\pdfoutput=1

\documentclass[11pt,letter,subeqn,fleqn]{article}
\usepackage{amsmath, amssymb,amsthm,cite}
\usepackage[normal]{subfigure}
\usepackage{subfig}
\usepackage{graphicx}
\usepackage{float}
\usepackage[utf8]{inputenc}
\usepackage{color}
\usepackage{url}

\usepackage{epstopdf}
\usepackage{enumitem}

\usepackage{array}
\newcolumntype{L}[1]{>{\raggedright\let\newline\\\arraybackslash\hspace{0pt}}m{#1}}
\newcolumntype{C}[1]{>{\centering\let\newline\\\arraybackslash\hspace{0pt}}m{#1}}
\newcolumntype{R}[1]{>{\raggedleft\let\newline\\\arraybackslash\hspace{0pt}}m{#1}}

\title{A Solution Set-Based Entropy Principle for constitutive modeling in Mechanics}

\author{ 
J. Heß \footnotemark[1],~~A. F. Cheviakov \footnotemark[2]\vspace{0.5cm}\\
\small $^{\rm a}$\emph{Department of Fluid Dynamics, Technical University of Darmstadt, Darmstadt, Germany}\vspace{0.2cm}\\
\small $^{\rm b}$\emph{Department of Mathematics and Statistics, University of Saskatchewan, Saskatoon, Canada}\vspace{0.2cm}\\ }

\def\beq{\begin{equation}}
\def\eeq{\end{equation}}
\def\barr{\begin{array}{ll}}
\def\earr{\end{array}}

\def\sg#1{{\rm #1}}
\def\const{\hbox{\rm const}}

\providecommand{\keywords}[1]{\textbf{\textit{Keywords:}} #1}

{\theoremstyle{definition} 

\newtheorem{remark}{Remark}
}

\newcounter{tabnum}

\newtheorem{lemma}{Lemma}

\begin{document}

\title{A Solution Set-Based Entropy Principle for Constitutive Modeling in Mechanics}

%\titlerunning{A Solution Set-Based Entropy Principle for Constitutive Modeling in Mechanics}        % if too long for running head
%
%\author{Julian Heß         \and
%       Alexei F. Cheviakov}
%
%\authorrunning{J. Heß, A. F. Cheviakov} % if too long for running head
%
%\institute{J. He{\ss} \at
%	Chair of Fluid Dynamics\\
%              Technical University of Darmstadt, Darmstadt, Germany\\
%              \email{hess@fdy.tu-darmstadt.de}      %  \\
%%             \emph{Present address:} of F. Author  %  if needed
% \and
%   Alexei F. Cheviakov \at
%Department of Mathematics and Statistics,\\
% University of Saskatchewan, Saskatoon, Canada \\
% \email{shevyakov@math.usask.ca}
%}


\maketitle

\begin{abstract}
Entropy principles based on thermodynamic consistency requirements are widely used for constitutive modeling in continuum mechanics, providing physical constraints on \emph{a priori} unknown constitutive functions. The well-known M{\"u}ller-Liu procedure is based on Liu's lemma for linear systems. While the M{\"u}ller-Liu algorithm works well for basic models with simple constitutive dependencies, it cannot take into account nonlinear relationships that exist between higher derivatives of the fields in the cases of more complex constitutive dependencies.

The current contribution presents a general solution set-based procedure, which, for a model system of differential equations, respects the geometry of the solution manifold, and yields a set of constraint equations on the unknown constitutive functions, which are necessary and sufficient conditions for the entropy production to stay nonnegative for any solution of the model. Similarly to the M{\"u}ller-Liu procedure, the solution set approach is algorithmic, its output being a set of constraint equations and a residual entropy inequality. The solution set method is applicable to virtually any physical model, allows for arbitrary initially postulated forms of the constitutive dependencies, and does not use artificial constructs like Lagrange multipliers. A \verb|Maple| implementation makes the solution set method computationally straightforward and useful for the constitutive modeling of complex systems.

%***Same as Liu if entropy ineq is linear in highest ders, and can also correctly handle cases when nonlinear, and/or when highest derivatives are related, by taking differential consequences of the field equations into account. ***

Several computational examples are considered, in particular, models of gas, anisotropic fluid, and granular flow dynamics. The resulting constitutive function forms are analyzed, and comparisons are provided. It is shown how the solution set entropy principle can yield classification problems, leading to several complementary sets of admissible constitutive functions; such classification problems have not previously appeared in the constitutive modeling literature.

\keywords{Constitutive modeling \and Continuum mechanics \and Entropy principle \and Solution set \and Symbolic computations}
\end{abstract}
%----------------------------------------------------------------------------------------------------------------------------------------------------------------------------------
\section{Introduction}

Models describing physical continua are commonly formulated as sets of equations which express balance laws and equations of state based on underlying physical principles. Such general models are often given systems of partial differential equations (PDE), containing \emph{free elements}, namely, constitutive functions and constant parameters. A specification of a particular form or class of such arbitrary elements corresponds to the description of a specific physical situation or a specific medium, thus a constriction to a certain material behavior.

The dependencies of such constitutive functions, called the \emph{constitutive dependencies} in the following, can be rather general, involving independent variables of the problem and physical fields, as well as their derivatives to a certain order. To describe a class of materials and mathematically analyze a class of models, it is often necessary to narrow down the constitutive dependencies, providing specific, practically relevant forms of their constitutive functions. This constitutes the main problem of \emph{constitutive modeling}. The principles for the derivation of constitutive function forms include, for example, the use of experimental data and of course theoretical considerations, such as physical conservation principles, geometric symmetry requirements \cite[p. 155ff]{hutter_continuum_2004}, the use of Lie symmetries \cite{magnenet2005lie,magnenet2009new,magnenet2014symmetry,magnenet2009continuous,rahouadj2003thermodynamic}, homogenization-based methods \cite{goda2012micropolar,goda20143d} or thermodynamical requirements \cite{coleman_thermodynamics_1963,kremer1986extended,liu_method_1972,muller_thermodynamic_1968,muller_thermodynamics_1984}.

The current work continues a series of papers that apply thermodynamic principles to the problem of constitutive modeling in continuum mechanics. This approach was pioneered by Ingo M{\"u}ller \cite{muller_thermodynamic_1968,muller1967entropy}, who, picking up the ideas of Coleman and Noll \cite{coleman_thermodynamics_1963} in the framework of Truesdell's thermodynamics \cite{truesdell1957sulle,truesdell_mechanical_1962}, formulated a generalized entropy inequality for continuum models of mixtures. In particular, he related energy and entropy supply terms, thus coupling the PDEs of the model to the basic entropy inequality, and introduced an extended form of the entropy inequality, which involved parts of the energy and momentum equations. Based on the requirement that the local entropy production must be nonnegative for all solutions, it followed that coefficients at ``free'', highest-order partial derivatives of physical fields not present in the constitutive dependencies must vanish. This led to \emph{constraint equations} restricting possible forms and limiting the dependencies of constitutive functions of the problem.

While M{\"u}ller's initial framework was restricted to rather special problems in physics of mixtures, a generalization was proposed by Liu \cite{liu_method_1972,liu_phd1972irreversible}. It is based on the Liu's Lemma (Section \ref{sec:MuellerLiuAppr}) holding for linear algebraic inequalities, and is related to the Farkas-Minkowski inequality lemma and the Hahn–Banach separation theorem \cite{hauser2002historical}. Liu's lemma  was then applied \cite{liu_method_1972} to nonlinear PDE models and entropy inequalities: For a given entropy inequality, the goal is to eliminate the highest derivatives of the fields, by ``subtracting off'' the model equations multiplied by so-called Lagrange multipliers, treated as additional constitutive functions. This yields a generalized entropy inequality, which still has some derivatives not present in the constitutive dependencies. At the terms linear in such derivatives, the coefficients then are set to zero, which leads to a set of constraints -- the \emph{Liu identities}. The latter are interpreted as constraints on the system's constitutive functions. This, in short, is currently known as \emph{the classical M{\"u}ller-Liu procedure.} This approach has been used for constitutive modeling in many applications, including chemical processes \cite{reis2016two}, viscoelastic bodies \cite{liu2008entropy}, materials in electromagnetic fields \cite{liu1972thermodynamicsferro,svendsen2005continuum} and granular flows \cite{hess2017thermodynamically,schneider2009solid}. The M{\"u}ller-Liu procedure is fully algorithmic, and has recently been implemented in symbolic software \cite{cheviakov2017symbolic}.

From the mathematical point of view, in order to study the properties of solutions to nonlinear PDEs rather than solutions of linear algebraic equations, one must take into account the structure of their solution manifold in a jet space of variables, physical fields, and their derivatives. Indeed, there is a substantial difference between linear equations or inequalities, where solutions and isosurfaces are given by linear, affine spaces, and inequalities that must hold on the solution manifolds of nonlinear PDEs in the jet spaces (see, e.g., \cite{olver2000applications}). The solution set of a PDE system in the appropriate jet space is described by the relationships between the field variables and their derivatives, given by the original dynamic equations and their \emph{differential consequences}. The latter are never computed in the M{\"u}ller-Liu algorithm. As a result, for generic constitutive dependencies, the Liu identities \emph{do not provide necessary conditions} for the local entropy production to be nonnegative, but may be overly restrictive, as is shown in examples below. It is also possible that Liu identities are not sufficient to guarantee nonnegative entropy production; this would happen, for example, if the postulated dependencies of Lagrange multipliers are overly general, and so they may be chosen to eliminate ``too much" of the entropy inequality, or perhaps the entropy inequality as a whole, without any restrictions on the forms of the actual constitutive functions of the problem.
%
%A simple explanation for that statement is that the M{\"u}ller-Liu algorithm sets to zero coefficients at different higher-order derivatives that may not be independent, and yields, figuratively speaking, $x=0, y=0$ instead of $x+y=0$.

In contrast to the M{\"u}ller-Liu algorithm, M{\"u}ller's original approach in the earlier works \cite{muller_thermodynamic_1968,muller1967entropy} is physical, and proceeds through the elimination of terms in the entropy inequality, proposing a relation to similar terms in the dynamic model PDEs \emph{on the solution set}. In \cite{muller1970new}, M{\"u}ller enhanced the method, suggesting a revised approach that couples the system's PDEs via the highest derivatives of the field variables. While this method apparently was only taken up by M{\"u}ller in two subsequent works \cite{muller1971kaltefunktion,muller_coldness_1971}, the approach proposed by Liu in Ref.~\cite{liu_method_1972}, i.e. the classical M{\"u}ller-Liu procedure, has become generally accepted.

The principal result of this paper is a re-formulation of the entropy principle in a physically and mathematically sound form, in the spirit of the approach suggested by M{\"u}ller in his later works \cite{muller1971kaltefunktion,muller1970new,muller_coldness_1971}. The proposed algorithm, called the \emph{solution set approach}, offers the technical flexibility of Liu's algorithm (splitting of the problem by setting to zero the coefficients at ``free'' higher-order derivatives), but provides a \emph{necessary and sufficient condition} for the entropy inequality to hold, through the consideration of the entropy inequality on the solution manifold of the given model. The latter is achieved through the use of substitutions of a set of \emph{leading derivatives} of the model PDEs and their differential consequences. In other words, the proposed algorithm is a mathematically justified nonlinear generalization of the classical M{\"u}ller-Liu procedure, related to M{\"u}ller's revised approach in his mostly neglected work \cite{muller1970new}, providing a more precise, and usually less restrictive (compared to the M{\"u}ller-Liu procedure), set of constraints on the constitutive functions of interest. It should be noted that the solution set approach does not require the use of the artificial Lagrange multipliers in any form, thus the number of unknown functions compared to the M{\"u}ller-Liu algorithm is significantly reduced. Moreover, the simplification of constraint equations following from the solution set approach allows for \emph{classification} (case splitting) that leads to separate families of admissible constitutive dependencies.

Similarly to the original M{\"u}ller-Liu algorithm, the newly proposed solution set-based approach is fully algorithmic, and is implemented using the existing \verb|GeM| package for \verb|Maple| computer algebra system \cite{cheviakov2007gem,cheviakov2010computation,cheviakov2010symbolic}. The \verb|GeM| package, originally created for symbolic calculations of Lie-type symmetries and conservation laws of differential equations, represents a PDE system in a symbolic form, which allows for efficient collection of coefficients in (differential) polynomial expressions, and thus the automated generation of constraint equations on the constitutive functions. The constraint equations can be consequently simplified using the efficient \verb|rifsimp| routine \cite{reid1996reduction}, and possibly solved explicitly using a built-in \verb|Maple| PDE solver.

The paper is organized as follows. In Section \ref{sec:EntropyPrinciple:ML}, we briefly review the constitutive modeling principles, containing M{\"u}ller's original approach and the M{\"u}ller-Liu algorithm, before we introduce the solution set-based entropy principle (Section \ref{sec:Entr:GenD}), using the model of a simple heat-conducting compressible anisotropic fluid as a running example.

In Section \ref{sec:Maple:EG1}, a \verb|Maple| implementation of the solution set method is presented and illustrated by an elementary example, taken from one-dimensional gas dynamics. It is shown how a classification problem arises, leading to several complementary classes of admissible constitutive functions. The results are compared with those for the classical M{\"u}ller-Liu procedure.

In Section \ref{sec:eg2}, a second, more involved computational example is considered: a model of a two-dimensional heat-conducting fluid. We show, for a simple choice of constitutive dependencies, the same constraint equations for the solution set approach and the classical M{\"u}ller-Liu procedure are derived. The solution set approach, however, can be used to take the problem further, and eliminate the residual entropy inequality, leading to a set of cases that correspond to the locally adiabatic motion of the fluid (Section \ref{sec:eg2:adiab}).

The third example considers a non-simple two-dimensional heat-conducting fluid, where the constitutive dependence involves a time derivative of a physical field (Section \ref{sec:Maple:EG3}). We show how the M\"{u}ller-Liu procedure leads to an overdetermined, overly restrictive set of Liu identities, while the solution set algorithm behaves consistently, yielding a set of necessary and sufficient conditions for the nonnegative entropy production. The paper is concluded with a discussion in Section \ref{sec:discussion}.

Appendix \ref{sec:Maple:EG4} discusses a more demanding computational example: a two-dimensional anisotropic granular solid flow model \cite{wang_comparison_1999,wang_shearing_1999}. It is shown how constraint equations and the residual entropy inequality are obtained. In this example, differential consequences of the model equations are considered as well as symmetrization conditions, required by the selected constitutive dependence.

%	
%-------------------------------------------------------------------------------------------------------------------------------------------------
\section{Constitutive modeling and the entropy principles}\label{sec:EntropyPrinciple:ML}

Here, we give a short account of the problem of constitutive modeling, in particular of the approaches based on thermodynamic entropy principles. For more details, see \cite{cheviakov2017symbolic} and references therein, as well as \cite{hutter_continuum_2004}.

% and the approaches proposed by M{\"u}ller and Liu, see \cite{muller_thermodynamic_1968}, \cite{liu_method_1972,muller_thermodynamics_1984} and \cite{muller1970new}. For a more detailed discussion, we refer to \cite{cheviakov2017symbolic}.

\subsection{Constitutive modeling, entropy inequality, and the running example}\label{sec:ConstitModProb}

Consider the motion of a physical continuum within a domain $\Omega\in \mathbb{R}^n$, $n\geq 1$, in time $t$, with spatial coordinates $x = (x_1,\ldots, x_n)\in \mathbb{R}^n$. The dependent variables are the physical parameters of the system, stated as fields
\[
\phi=(\phi_1(t, x),\ldots,\phi_m(t, x)).
\]
The evolution of these fields is assumed to be described by a system of equations $\Pi^{\phi}$, consisting generally of algebraic, ordinary and/or partial differential equations
\[
\Pi^{\phi} = \left\{\Pi^{\phi_1}, \Pi^{\phi_2},\ldots,\Pi^{\phi_m}\right\},
\]
and possible additional algebraic or differential constraints.
In addition to the fields, the model equations contain one or more \emph{constitutive functions}
\beq\label{eq:const:dep}
\psi=\psi(\phi_C).
\eeq
The \emph{constitutive dependencies} $\phi_C$ determine the choice of relevant arguments of the constitutive functions; they may include dependent and independent variables of the problem, as well as first and higher order derivatives of the dependent variables. \footnote{Processes with memory, in which constitutive dependencies may include integral- and delay-based quantities, and other complex processes requiring the use of nonlocal terms, are not considered in this contribution.}

%The specification of particular classes of dependencies \eqref{eq:const:dep} relevant to specific applications is the essence of the \emph{problem of constitutive modeling}.

While the main set of balance equations $\Pi^{\phi}$ may describe the behavior of a broad range of materials, the specification of the constitutive functions
%provides the closure of the system and
defines the material behavior of a given medium. The problem of constitutive modeling consists in the formulation of physically relevant constitutive dependencies and the description of how they are incorporated in the constitutive functions. The latter is often formulated in terms of \emph{constraint equations}. The principles of constitutive modeling can include general assumptions like invariance with respect to coordinate transformations (the principle of material objectivity), geometrical symmetries, material isotropy, as well as other physical postulates and simplifying assumptions of physical and mathematical nature, see, e.g. \cite{cheviakov2017symbolic,hutter_continuum_2004} and the references therein.

As a running example for the purpose of depiction and comparison, we will use the model of a simple heat-conducting compressible anisotropic fluid. The physical fields here are given by the density $\rho=\rho(t,x)$, the velocity $v_i=(v_1(t,x),\ldots,v_n(t,x))$, and the temperature $\vartheta=\vartheta(t,x)$; $x=(x_1,\ldots,x_n)\in \mathbb{R}^n$ is the spatial variable. The constitutive functions of the model are the fluid stress tensor $T_{ij}=T_{ij}(t,x)$, the outgoing energy flux vector $q_i=(q_1(t,x),\ldots,q_n(t,x))$, and the internal energy per unit mass $\epsilon=\epsilon(t,x)$. The equations
\begin{subequations}\label{eq:fluid:gen}
\begin{equation}\label{eq:fluid:gen:dens}
\Pi^{\rho}:~ \partial_t \rho + \partial_i \left(\rho v_i \right)=0,
%\Pi^{v}: \rho D_t v - (\div T) =0,\qquad i=1\ldots, n,\\
\end{equation}
\begin{equation}\label{eq:fluid:gen:mom}
\Pi_i^{v}:~\rho\, \sg{D}_t v_i - \partial_j  T_{ij}- \rho g_i =0,\qquad i=1\ldots, n,\\[2ex]
\end{equation}
\begin{equation}\label{eq:fluid:gen:en}
\Pi^{\epsilon}:~\rho \,\sg{D}_t \epsilon +  \partial_i q_i -  T_{ij} \partial_j  v_i =\rho r,
\end{equation}
\end{subequations}
are the balance of mass, of momentum and (internal) energy, respectively. For this example, the constitutive functions $\psi=\{T_{ij}, \epsilon, q_i\}$  may depend, for example, on $\phi_C =  \left\{\rho, \vartheta\right\}$, see equation \eqref{eq:const:dep}. As a main postulate for constitutive modeling based on entropy principles, the local form of the \emph{entropy inequality} is introduced as
\begin{equation}\label{eq:entropy_ineq:Liufluid}
\Pi^{\eta}:~\rho \,\sg{D}_t\, \eta + \partial_i \Phi_i  - \rho s \geq 0,
\end{equation}
and assumed to hold for all physical solutions of the model based on the dynamic equations \eqref{eq:fluid:gen}. The outgoing entropy flux vector $\Phi = (\Phi_1(t,x),\ldots,\Phi_n(t,x))$ and the entropy per unit mass $\eta=\eta(t,x)$ are additional constitutive functions introduced herewith.

In Eqs. \eqref{eq:fluid:gen}-\eqref{eq:entropy_ineq:Liufluid}, $\rho r$ denotes the external energy supply density, $\rho g_i$ is the body force (gravity), and $\rho s$ is the entropy supply. These external supplies are considered to be independent of the internal material behavior, and thus are not constitutive functions. For the sake of brevity, we will denote time and space partial derivatives as
\[
\frac{\partial }{\partial t}=\partial_t,\qquad \frac{\partial}{\partial x_i}=\partial_i, %\qquad \frac{\partial^2}{\partial x_i \partial x_j}=\partial^2_{i,j},
\]
%(***latter never used!***)
and use the Einstein summation convention where appropriate. The material derivative operator used in Eqs. \eqref{eq:fluid:gen}-\eqref{eq:entropy_ineq:Liufluid} is given by
\beq \label{eq:TotD}
\sg{D}_t = \frac{\partial}{\partial t} +v_i \frac{\partial}{\partial x_i}= \partial_t +v_i \partial_i.
\eeq
Where appropriate, we will use the Cartesian coordinate names $x\equiv x_1$, $y\equiv x_2$, $\ldots$, and the respective velocity components $u\equiv v_1$, $v\equiv v_2$ and $w\equiv v_3$.

\begin{remark}
Consider the global entropy production rate $Q(t)$ in a simple bounded domain $V(t)$ with a piecewise-smooth boundary $\partial V(t)$, moving with the medium:
\beq \label{eq:glob:Q}
Q_V(t) = \frac{d}{dt} \int_{V(t)} \rho \eta \, dV +\oint_{\partial V(t)} \Phi_i dS_i - \int_{V(t)} \rho s \, dV.
\eeq
Note that this construct works the same whether $V(t)$ is a volume in 3D or an area in 2D; here, $\Phi_i dS_i = \Phi\cdot dS$, where $dS$ is directed along the outward normal of the domain boundary $\partial V(t)$. With the Leibniz's rule for moving domains, Eq. \eqref{eq:glob:Q} can be rewritten as
\beq \label{eq:glob:Q2}
Q_V(t) = \int_{V(t)} \left( \partial_t(\rho \eta) +  \partial_i (\rho \eta v_i + \Phi_i) -  \rho s \right)\, dV.
\eeq
The quantity in the brackets is the local entropy inequality plus the continuity equation \eqref{eq:fluid:gen:dens} times $\eta$, and hence the requirement \eqref{eq:entropy_ineq:Liufluid} is equivalent to the requirement
\[
\partial_t(\rho \eta) + \partial_i (\rho \eta v_i + \Phi_i)  - \rho s \geq 0.
\]
Consequently, the nonnegative local entropy production \eqref{eq:entropy_ineq:Liufluid} guarantees that the global entropy production \eqref{eq:glob:Q} is nonnegative, $Q_V(t)\geq 0$, in every moving material domain $V(t)$.

If the local entropy production \eqref{eq:entropy_ineq:Liufluid} vanishes identically, it means that the entropy rate of change $\frac{d}{dt} \int_{V(t)} \rho \eta \, dV$ in every material volume ${V(t)}$ is only due to the entropy fluxes and possibly sources (when $\rho s\ne 0$). The process is globally adiabatic when
\[
\frac{d}{dt} \int_{W} \rho \eta \, dV =0
\]
for the static problem domain $W$. In the absence of the local entropy sources, $s=0$, a globally adiabatic process is may be described by vanishing entropy inequality
\begin{equation}\label{eq:entropy_ineq:Liufluid:glob_adiab}
\rho \,\sg{D}_t\, \eta + \partial_i \Phi_i =0
\end{equation}
and an additional requirement that entropy fluxes and velocities vanish, $v_i=\Phi_i=0$, on the problem domain boundary $\partial W$.

\end{remark}

Below, we briefly describe M{\"u}ller's original constitutive modeling approach (Section \ref{sec:MuellerAppr}), the famous M{\"u}ller-Liu algorithm (Section \ref{sec:MuellerLiuAppr}) involving Lagrange multipliers $\Lambda^{\phi}(\phi_C)$, and the proposed solution set entropy principle (Section \ref{sec:Entr:GenD}). A comparative summary of these approaches is presented in Table \ref{tab:principles}.

\begin{table}[H]
  \centering
  \begin{tabular}{|C{5cm}||C{3cm}|C{3cm}|C{3cm}|}
     \hline
     \textbf{Algorithm} & \textbf{M{\"u}ller} & \textbf{M{\"u}ller-Liu} & \textbf{Solution set} \\ \hline\hline
     Applicability & Continuum mechanics  & Any model & Any model \\ \hline
     Constitutive functions & $\psi(\phi_C)$ & $\psi(\phi_C)$, $\Lambda^{\phi}(\phi_C)$ & $\psi(\phi_C)$ \\ \hline
     Takes into account differential consequences? & No & No  &  Yes\\  \hline
     Ensures nonnegative entropy production for all solutions? & Yes, possibly with a residual entropy inequality  & Not generally  & Yes, possibly with residual entropy inequality  \\ \hline
     Residual entropy inequality  & Present or absent & Present or absent, involves Lagrange multipliers $\Lambda^{\phi}$ &  Present or absent\\  \hline
   \end{tabular}
  \caption{A comparison of M{\"u}ller's, the M{\"u}ller-Liu, and the solution set entropy principle.}\label{tab:principles}
\end{table}

%
%----------------------------------------------------------------------------------------------
\subsection{M{\"u}ller's constitutive modeling approach }\label{sec:MuellerAppr}

M{\"u}ller \cite{muller_thermodynamic_1968}, following and generalizing the work of Coleman and Noll \cite{coleman_thermodynamics_1963}, argued that the entropy inequality \eqref{eq:entropy_ineq:Liufluid} has to be regarded as a restriction on the constitutive functions rather than as a restriction on the processes. This means that, for the inequality to hold for every possible process, for example for every physical solution $\rho,v_i,\vartheta$ of Eqs. \eqref{eq:fluid:gen}, adequate constitutive functions have to be chosen. In \cite{muller1970new}, M{\"u}ller states the general framework of his approach\footnote{\label{foot:1}We note that while M{\"u}ller, in his original approach \cite{muller_thermodynamic_1968}, considered a mixture consisting of several constituents, we use simpler examples of single-phase materials to outline the ideas.} as
%\begin{lemma}[M{\"u}ller] \label{lem:mueller}
\begin{enumerate}[label=(\roman*)]
\item A process is defined as a set of functions $\rho(t,x), \chi(t,x), \vartheta(t,x)$ representing
the density, the motion (instead of velocity) and the (empirical) temperature.
\item Constitutive equations are formulated for the stress tensor $T_{ij}$, the internal
energy $\epsilon$ and the heat flux $q_i$, such that a set $T_{ij},\epsilon,q_i$ corresponds to every process.
\item A thermodynamic process is defined as a process which is a solution of the
equations of balance for the mass, the momentum and the energy.
\item In a body $\mathcal{B}$, there exists a scalar extensive quantity
which cannot decrease in any thermodynamic process, if its flux through
the surface of $\mathcal{B}$ vanishes, and whose density $\eta$ and flux $\Phi_i$ are determined by
constitutive relations. This quantity is called entropy. (Entropy Principle)
\end{enumerate}

M{\"u}ller's principle states that all physically relevant processes for a system of PDEs as in \eqref{eq:fluid:gen} must satisfy the entropy inequality, see Eq. \eqref{eq:entropy_ineq:Liufluid}. M{\"u}ller then proceeds to relate the entropy supply term with the external radiative energy supply, so that $s = {r_R}/{\vartheta}$, coupling the energy balance equation with the entropy inequality. In a similar way, through a term accounting for the body force density, the momentum balance equations are coupled to the energy balance and, in turn, incorporated into the entropy inequality. These substitutions lead to an \emph{extended} entropy inequality, a linear combination that involves essential parts of the energy and the momentum balance equations. With several assumptions on the form of the constitutive functions and their dependence on the fields and the derivatives of the fields, the constitutive functions are substituted in the extended entropy inequality. The resulting inequality is linear in the highest derivatives of the field variables, which can take any value for various solutions. To ensure nonnegative entropy production, the corresponding coefficients of these highest derivatives must consequently vanish. This gives a set of restrictions in the form of PDEs, providing constraints on the constitutive functions.

The original form of M{\"u}ller's approach outlined above is mathematically and physically sound, yet applicable only to models described by PDEs with source and forcing terms of rather special structure, and thus is generally not applicable to a broad range of constitutive modeling problems (see Table \ref{tab:principles}).

%
%
%----------------------------------------------------------------------------------------------
\subsection{Liu's approach. The classical and generalized M{\"u}ller-Liu algorithm}\label{sec:MuellerLiuAppr} %Principle

Liu \cite{liu_method_1972,liu_phd1972irreversible,muller_thermodynamics_1984} generalized M{\"u}ller's approach via the consideration of Lagrange multipliers and a constrained entropy inequality, proposing what is known as the M{\"u}ller-Liu entropy principle. It is based on the following lemma stated for linear algebraic equations.
\begin{lemma}[Liu] \label{lem:Liu}
Let $z\in \mathbb{R}^p$, and let $M$ be a $p\times n$ real matrix. Consider a linear system $MY +z = 0$ of $p$ equations on the components of the unknown vector $Y\in \mathbb{R}^n$, with a non-empty solution set $S$. Let also $\mu\in \mathbb{R}^n$, $\mu\ne 0$, and $\zeta\in \mathbb{R}$ be given. Then the following statements are equivalent:
\begin{enumerate}
  \item $\forall Y\in S$, $\mu^T Y + \zeta \geq 0$;
  \item $\exists \lambda \in \mathbb{R}^p$ such that $\forall Y \in \mathbb{R}^n$, $\mu^T Y + \zeta  - \lambda^T(MY +z) \geq 0$;
  \item $\exists \lambda \in \mathbb{R}^p$ such that $\mu=M^T\lambda$, and $\zeta  \geq \lambda^T z$.
\end{enumerate}
\end{lemma}

In contrast to M{\"u}ller's procedure \cite{muller_thermodynamic_1968}, in Liu's approach, all external source terms are neglected, and the coupling of the set of balance equations is achieved through the introduction of so-called Lagrange multipliers (even though no actual Lagrangian is considered). Instead of a sequence of substitutions, every equation $\Pi^{\phi}$ of the given system, e.g. Eqs. \eqref{eq:fluid:gen}, is regarded as a constraint on the entropy inequality, with a respective Lagrange multiplier $\Lambda^{\phi}$. A linear combination
\beq\label{eq:ML:extEntr}
\widetilde{\Pi}^{\eta}:~\Pi^{\eta} - \Lambda^{\phi} \Pi^{\phi} \geq 0
\eeq
is called the \emph{extended entropy inequality}. Again, the constitutive functions with their postulated dependencies $\phi_C$ are substituted into Eq. \eqref{eq:ML:extEntr}, and the coefficients of the highest derivatives are set to zero, in accordance with Lemma \ref{lem:Liu}. This yields constraint equations, called the \emph{Liu identities}, and a residual entropy inequality, as constraints on the constitutive functions. This entropy principle in its formulation of M{\"u}ller and Liu has been used in numerous works to the present, including both articles, e.g. \cite{hess2017thermodynamically,liu2008entropy,reis2016two,svendsen_thermodynamics_1995} as well as extensive books \cite{hutter_continuum_2004,liu2002shih,muller1998rational,schneider2009solid}.

According to equation \eqref{eq:ML:extEntr}, the M{\"u}ller-Liu approach is based on the linear superposition of the entropy inequality and the governing equations. It is appropriate only in the special case of linearly occurring highest derivatives, and there is no relationship between such derivatives. With this, the linear algebra-based Lemma \ref{lem:Liu} holds for the differential equations, and the Liu identities yield meaningful restrictions on constitutive functions. This is the case, for instance, in our computational examples 1 and 2 below, i.e. in Sections \ref{sec:Maple:EG1} and \ref{sec:eg2}.

In general, however, the solution set of the physical equations in the jet space of independent variables, dependent variables, and derivatives of dependent variables, is not a linear variety, and Lemma \ref{lem:Liu} is not applicable. Instead of trying to select Lagrange multipliers that cancel ``free'' derivatives, one must take into account the actual relationships between fields and their derivatives in the jet space. Such relationships are given by the equations themselves and their differential consequences, as explained in Section \ref{sec:Entr:GenD} below (see also Table \ref{tab:principles}). As a result, for a general physical model, the M{\"u}ller-Liu approach yields constraint equations that are neither necessary nor sufficient for a nonnegative entropy production -- see the  computational example 3  in Section \ref{sec:Maple:EG3} below.

It is worth noting that within the M{\"u}ller-Liu approach, it is not clear how to choose the dependence of the Lagrange multipliers $\Lambda^{\phi}$. The heuristic rule in the literature is to postulate $\Lambda^{\phi}= \Lambda^{\phi}(\phi_C)$, the same constitutive dependence as for the unknown constitutive functions; it is not clear whether this is necessary, and what would be sufficient. In fact, if the set $\phi_C$ is broad enough, one might choose, for example, a singular multiplier $\Lambda^{\phi}=\Pi^{\eta}/\Pi^{\rho}$, and thus completely annihilate the extended entropy inequality: $\widetilde{\Pi}^{\eta}\equiv 0$. In this or similar cases, the resulting Liu identities would provide insufficient conditions for nonnegative entropy production.

\smallskip

An extension of the M{\"u}ller-Liu procedure, and a corresponding symbolic computer algebra implementation, were presented in Ref.~\cite{cheviakov2017symbolic}, where, in the extended entropy inequality \eqref{eq:ML:extEntr}, the coefficients were set to zero not only at the terms linear in a set of highest derivatives, but also at general polynomial involving other derivatives that did not enter the constitutive dependency $\phi_C$.

%-------------------------------------------------------------------------------------------------------------------------------------------------
\subsection{The solution set entropy principle}\label{sec:Entr:GenD}
%Generalized

As it can be seen from the comparison of M{\"u}ller's and Liu's approaches above, Liu's algorithm employs essentially the same idea of using dynamic equations to compose an extended inequality, annihilating certain terms in the entropy inequality. Then, Liu's algorithm proceeds by considering a certain set of higher derivatives ``free'', and thus, through Lemma \ref{lem:Liu}, sets the respective coefficients to zero, obtaining a set of restrictions on the constitutive quantities and a residual entropy inequality.

While Liu's Lemma itself is certainly correct, its application to describe relationships between a nonlinear inequality constraint and a solution manifold of a nonlinear PDE model in the corresponding jet space is not mathematically justified. There exists a significant difference between linear equations/inequalities, with solutions and isosurfaces being represented by linear or affine spaces, and solution manifolds of nonlinear differential equations. Here, it is normally assumed that a given PDE system is \emph{locally solvable} (e.g., \cite{olver2000applications}), that is, the solution set of the given PDE system in the jet space is actually represented by these PDEs.

In addition to the model PDEs, for the example governed by Eqs. \eqref{eq:fluid:gen}, further relationships between the field variables are provided by their \emph{differential consequences}. It follows that in the M{\"u}ller-Liu idea of eliminating the highest derivatives of the fields, the extended entropy inequality $\widetilde{\Pi}^{\eta}$ of Liu's algorithm may need to involve the differential consequences of the equations $\Pi^{\phi}$, with additional multipliers. For example, for the heat-conducting fluid model \eqref{eq:fluid:gen}, if one extends the constitutive dependence $\phi_C$ to include $\partial_t \rho$, it is clear that the time derivative $\partial_t \epsilon$ in the energy balance will involve the second derivative $\partial_t^2 \rho$. It follows that the entropy condition $\widetilde{\Pi}^{\eta}$ should include an additional term $\Lambda^{\rho}_1 \, \partial_t \Pi^{\rho}$ involving a time differential consequence of the continuity equation. In a similar way, spatial differential consequences $\partial_{x_i} \Pi^{\rho}$ would be added, with appropriate Lagrange multipliers. Little-known but proposing the basis for his famous works on the coldness-function \cite{muller1971kaltefunktion,muller_coldness_1971}, M{\"u}ller suggested a revised approach that avoids these drawbacks and, instead of applying the method of Lagrange multipliers, solves the system's balance laws for the highest derivatives of the fields and inserts these in the entropy inequality in \cite{muller1970new}.

Below, we propose a generalized entropy principle algorithm, which proceeds by working directly on the solution set of the model in the manner of M{\"u}ller's revised approach in \cite{muller1971kaltefunktion,muller1970new,muller_coldness_1971}, at the same time retaining the general and systematic structure of Liu's algorithm. The principle uses the computation of a set of \emph{leading derivatives} and their differential consequences on the solution space, and their substitution into the entropy inequality, thus requiring no Lagrange multipliers. The same running example of Section \ref{sec:ConstitModProb}, in two spatial dimensions is used as an illustration.

As outlined in Table \ref{tab:principles}, the solution set approach supersedes the M{\"u}ller-Liu algorithm, providing a set of sufficient conditions to ensure that the entropy production is nonnegative on any solution of a model, for any posed constitutive dependences.

\subsection*{\bf The solution set algorithm and its illustration} %\label{sec:Entr:GenD1}

%***** The example is based on file \verb|Liu_fluid_dynamics_book_reduced(7.5)_SolvApproach5t.mw| ***

\begin{enumerate}
  \item For a given physical model, define the fields of interest $\phi$, and the dynamic PDEs $\Pi^{\phi}$. In the running example, we have the dependent variables
\[
\phi = (\rho, v, \vartheta),
\]
and the governing equations \eqref{eq:fluid:gen} with vanishing source terms, so that $g_i, r=0$\footnote{\label{foot:2}While Liu's assumption here is that the source terms cancel out each other, we follow the argumentation that the external influences do not affect the material constraints and thus can be neglected.}.

%--------------------------------------- 2
  \item Define the entropy inequality $\Pi^{\eta}$. For our example, it is given by \eqref{eq:entropy_ineq:Liufluid} with vanishing source terms ($s=0$).

%--------------------------------------- 3
  \item Define the constitutive functions $\psi$ of the model and postulate their dependencies \footnote{\label{foot:3}Note that constitutive dependencies for different constitutive functions may be different, unless the principle of equipresence is enforced, see \cite[p. 289f]{hutter_continuum_2004}}.

  In the running example, we choose the unknown constitutive functions to be
\beq\label{eq:simplefluid:constitF:our:algo}
\barr
    \psi=(\epsilon, q_i, T_{ij}, \eta, \Phi_{i}),\quad i,j=1,2,
\earr
\eeq
with dependencies
\beq\label{eq:simplefluid:constit:our:algo}
\barr
\epsilon=\epsilon(\rho,\vartheta),\quad \eta=\eta(\rho,\vartheta),\\[2ex]
T_{ij}=T_{ij}(\rho, \vartheta, \partial_i \vartheta), \quad q_i=q_i(\rho, \vartheta, \partial_i \vartheta), \quad \Phi_{i}=\Phi_{i}(\rho, \vartheta, \partial_i \vartheta).
\earr
\eeq
%  \beq\label{eq:simplefluid:constit:our:algo}
%  \barr
%   \epsilon=\epsilon(\rho,\vartheta),& \eta=\eta(\rho,\vartheta),\\
%   T_{ij}=T_{ij}(\rho, \vartheta, \partial_i \vartheta), & q_i=q_i(\rho, \vartheta, \partial_i \vartheta), \\ \Phi_{i}=\Phi_{i}(\rho, \vartheta, \partial_i \vartheta).
%  \earr
%  \eeq

%--------------------------------------- 4
  \item Substitute the chosen forms of the constitutive functions into the governing equations and the entropy inequality; carry out partial differentiations. Thus, obtain the governing equation forms $\widehat{\Pi}^{\phi}$ and the entropy inequality form $\widehat{\Pi}^{\eta}$.

  In the running example, for instance, the continuity equation \eqref{eq:fluid:gen:dens} has no constitutive functions and therefore does not change: $\widehat{\Pi}^{\rho} = \Pi^{\rho}$; the momentum equations \eqref{eq:fluid:gen:mom} become
  \[
  \widehat{\Pi}_i^{v}:~\rho \left(\partial_t v_i + v_j \,\partial_{j}v_i\right)  - \left(\dfrac{\partial  T_{ij}}{\partial \rho}\, \partial_{j}\rho  + \dfrac{\partial  T_{ij}}{\partial \vartheta}\, \partial_{j}\vartheta+ \dfrac{\partial  T_{ij}}{\partial (\partial_k \vartheta)}\, \partial_j \partial_k\vartheta\right) =0,\qquad i=1,2,
  \]
  the energy equation is given by
  \[
  \barr
  \widehat{\Pi}^{\epsilon}:~&\rho \left( \dfrac{\partial  \epsilon}{\partial \rho}\, \partial_t \rho + \dfrac{\partial  \epsilon}{\partial \vartheta} \,\partial_{t}\vartheta  + v_i \left[ \dfrac{\partial  \epsilon}{\partial \rho} \,\partial_{i}\rho + \dfrac{\partial  \epsilon}{\partial \vartheta}\, \partial_{i}\vartheta \right]\right)\\[3ex]
  &+\left(\dfrac{\partial  q_{i}}{\partial \rho}\, \partial_{i}\rho  + \dfrac{\partial  q_{i}}{\partial \vartheta}\, \partial_{i}\vartheta+ \dfrac{\partial  q_{i}}{\partial (\partial_k \vartheta)}\, \partial_i \partial_k\vartheta\right) -  T_{ij} \partial_{j} v_i  =0,
  \earr
  \]
   and the entropy inequality after the substitution becomes
  \beq\label{eq:running:Pi:Eta}
\barr
  \widehat{\Pi}^{\eta}:~&\rho \left( \dfrac{\partial  \eta}{\partial \rho} \,\partial_{t}\rho + \dfrac{\partial  \eta}{\partial \vartheta}  \,\partial_{t}\vartheta  + v_i \left[ \dfrac{\partial  \eta}{\partial \rho} \partial_{i}\rho +  \dfrac{\partial  \eta}{\partial \vartheta} \partial_{i}\vartheta \right]\right)\\[3ex]
&+ \left(\dfrac{\partial  \Phi_{i}}{\partial \rho}\, \partial_{i}\rho  + \dfrac{\partial  \Phi_{i}}{\partial \vartheta}\, \partial_{i}\vartheta+ \dfrac{\partial  \Phi_{i}}{\partial (\partial_k \vartheta)}\, \partial_i \partial_k\vartheta\right) \geq 0.
  \earr
  \eeq
  In the above formulas, we use $\partial_{i}\rho$, etc., to denote partial derivatives of field variables, and full partial derivative symbols for the partial derivatives of constitutive functions with respect to fields.

%--------------------------------------- 5
  \item For the PDEs $\widehat{\Pi}^{\phi}$ with substituted constitutive function forms, define a set of \emph{leading derivatives}
      \[
      \phi_L = \{D^{\alpha_i} \phi_i\}_{i=1}^m,
      \]
      with respect to which the PDEs $\widehat{\Pi}^{\phi}$ can be solved. Here, each $D^{\alpha_i}$ denotes a specific (possibly mixed) partial derivative of order $\alpha_i\geq 1$. Rewrite the system $\{\widehat{\Pi}^{\phi}\}$ in \emph{the solved form} $\{\widetilde{\Pi}^{\phi}\}$ with respect to these leading derivatives:
\beq\label{eq:gen:PhysPDEs:solvedform}
\widetilde{\Pi}^{\phi_i}:~~D^{\alpha_i} \phi_i = F_i,\quad i=1,\ldots, m;
\eeq
  in particular, the right-hand side $F_i$ of every equation in the solved form may be a function of $t,x,\phi$, and derivatives of $\phi$, but involves neither the leading derivatives $\phi_L$ nor their \textit{differential consequences}.

  Assuming that the given PDE system is \emph{locally solvable} in some domain in the jet space, see, e.g. \cite{olver2000applications}, it follows that an arbitrary specification of $t,x$, the field variables $\phi$, and their derivatives contained in $F_i$, at any given point of the system solution manifold in the jet space, leads to an unique specification of the values of the leading derivatives $\phi_L$, and the corresponding solution in fact exists.

  In the running example, one may choose, for example, the leading derivatives to be highest $x$-derivatives of the dependent variables: $\phi_L = \{ \partial_{1}\rho, \partial_{1} v_i, \partial_{1}^2\vartheta\}$, or the highest $t$-derivatives, $\phi_L = \{ \partial_t \rho, \partial_t v_i, \partial_t \vartheta\}$. Any choice of a set of leading derivatives would lead to equivalent results; each one of the two specific choices outlined above yields a \emph{Kovalevskaya form} of the given PDE system, see \cite{olver2000applications} for details. We proceed with the second choice, which leads to a shorter set of equations $\{\widetilde{\Pi}^{\phi}\}$:
\beq\label{eq:fluid:solved:specT}
\barr
\widetilde{\Pi}^{\rho}:~ \partial_t \rho = -\partial_j (\rho v_j);\\[2ex]

\widetilde{\Pi}^{v}_i:~ \partial_t v_i = -v_j\, \partial_{j}v_i + \dfrac{1}{\rho}\left(\dfrac{\partial  T_{ij}}{\partial \rho}\, \partial_{j}\rho  + \dfrac{\partial  T_{ij}}{\partial \vartheta}\, \partial_{j}\vartheta+ \dfrac{\partial  T_{ij}}{\partial (\partial_k \vartheta)}\, \partial_j \partial_k\vartheta\right),~~~i=1,2;\\[2ex]

  \widetilde{\Pi}^{\epsilon}:~ \partial_{t}\vartheta = \dfrac{1}{\rho \, \partial  \epsilon/ \partial \vartheta}\left( - \dfrac{\partial  \epsilon}{\partial \vartheta}\,\rho v_i\,\partial_{i}\vartheta + \dfrac{\partial  \epsilon}{\partial \rho} \rho^2 \partial_j v_j\right.\\[2ex]
  \qquad\qquad
  \left.
  -\left[\dfrac{\partial  q_{i}}{\partial \rho}\, \partial_{i}\rho  + \dfrac{\partial  q_{i}}{\partial \vartheta}\, \partial_{i}\vartheta+ \dfrac{\partial  q_{i}}{\partial (\partial_k \vartheta)}\, \partial_i \partial_k\vartheta\right] +  T_{ij} \partial_{j} v_i\right).
\earr
\eeq

%--------------------------------------- 6
  \item Substitute the equations in the solved form, $\{\widetilde{\Pi}^{\phi}\}$, into the entropy inequality $\widehat{\Pi}^{\eta}$, thus enforcing it to hold on solutions. The resulting
            \emph{entropy inequality on solutions}
  \beq\label{eq:gen:EntrOnSol}
  \widetilde{\Pi}^{\eta}:~Q(Z)\geq 0
  \eeq
  is a nonlinear function of the constitutive functions $\psi$ \eqref{eq:simplefluid:constitF:our:algo} and the
  \emph{free elements}
  \beq\label{eq:gen:FreeEl:Z}
  Z=\{t,x,\phi, \partial_t \phi, \partial_i \phi, \partial_t^2 \phi,\ldots\} \setminus \phi_L \setminus \phi_C,
  \eeq
a set of possibly all independent variables, physical fields, and the derivatives of the latter,
excluding the leading derivatives $\phi_L$ (and where applicable, all of their differential consequences), and also excluding the arguments of the constitutive functions $\phi_C$.
Since the elements of $Z$ can be varied independently, the entropy inequality $\widetilde{\Pi}^{\eta}$ can be \emph{split} with respect to the free elements: all the coefficients of all independent powers/combinations of the elements of $Z$ in the entropy inequality $\widetilde{\Pi}^{\eta}$ must vanish separately.

  In the heat-conducting fluid example, the entropy  inequality on solutions $\widetilde{\Pi}^{\eta}$ \eqref{eq:gen:EntrOnSol} can be split with respect to the 23 free elements given by
  \beq\label{eq:gen:FreeEl:Z1:eg}
  {Z}=\{t, x, v_i, \partial_j v_i, \partial_j \rho, %\partial_j \vartheta,
  \partial_j\partial_k v_i, \partial_j\partial_k \rho, \partial_j\partial_k \vartheta\}.
  \eeq

%--------------------------------------- 7
  \item In the entropy inequality on solutions $\widetilde{\Pi}^{\eta}$, collect coefficients at all different polynomial terms involving elements of $\widetilde{Z}$. Set them to zero. Obtain a set of \emph{constraint identities}, which replace Liu identities as constraint equations on constitutive function.

       The term in $\widetilde{\Pi}^{\eta}$ with no factors from $\widetilde{Z}$ defines a \emph{residual entropy inequality}.

       For the running example, one has eight constraint equations, given by
\beq\label{eq:OurAlgo:RunEx:ConstrID}
\barr
\dfrac{\partial \epsilon}{\partial\vartheta} \dfrac{\partial \Phi_i}{\partial\rho} -\dfrac{\partial \eta}{\partial\vartheta} \dfrac{\partial q_i}{\partial\rho}=0,\quad i=1,2;\\[3ex]

\rho^2\left(\dfrac{\partial\epsilon}{\partial\rho}\dfrac{\partial \eta}{\partial\vartheta} - \dfrac{\partial \epsilon}{\partial\vartheta} \dfrac{\partial \eta}{\partial\rho} \right)  \delta_{ij}+\dfrac{\partial \eta}{\partial\vartheta} T_{ij}=0,\quad i,j=1,2;\\[3ex]

\dfrac{\partial \epsilon}{\partial\vartheta} \left(\dfrac{\partial \Phi_i}{\partial (\partial_j\,\vartheta)}+\dfrac{\partial \Phi_j}{\partial (\partial_i\,\vartheta)}\right)-\dfrac{\partial \eta}{\partial\vartheta} \left(\dfrac{\partial q_i}{\partial (\partial_j\,\vartheta)}+\dfrac{\partial q_j}{\partial (\partial_i\,\vartheta)}\right)=0, \quad i,j=1,2,
\earr
\eeq
with an additional assumption $\partial  \epsilon/ \partial \vartheta\ne 0$ due to the solved form \eqref{eq:fluid:solved:specT}. The residual entropy inequality is given by
\beq\label{eq:Liu:RunEx:ResIneq}
\left(\dfrac{\partial \epsilon}{\partial\vartheta}\right)^{-1}\,\partial_i \vartheta\, \left(  \dfrac{\partial \epsilon}{\partial\vartheta} \dfrac{\partial \Phi_i}{\partial\vartheta} -\dfrac{\partial \eta}{\partial\vartheta} \dfrac{\partial q_i}{\partial\vartheta} \right) \geq 0.
\eeq
If we compare \eqref{eq:Liu:RunEx:ResIneq} with the residual entropy inequality for the same example in \cite{cheviakov2017symbolic}, tackled with the M{\"u}ller-Liu entropy principle, it becomes apparent that the results are identical for $\Lambda^{\epsilon} = \dfrac{\partial \eta}{\partial \vartheta} / \dfrac{\partial \epsilon}{\partial  \vartheta}$, which is one of the Liu identities derived there. Furthermore, all of the resulting Liu identities are reproduced, without Lagrange multipliers, in \eqref{eq:OurAlgo:RunEx:ConstrID}.

%--------------------------------------- 8
  \item The last step is to analyze and solve the constraint equations. Since the equations are underdetermined, and the solution depends on a number of arbitrary functions, it may not be always possible to obtain simple closed-form solutions of the constraint equations. Moreover, special cases may arise that yield particular families of constitutive functions that are admissible from the point of view of thermodynamic consistency. Please note that here, we present only most general results; for details, particular cases, and discussion, we refer to Section \ref{sec:eg2:comp} below. In the current example, for the most general case, there is no restriction on the forms of energy density $\epsilon=\epsilon(\rho,\vartheta)$ and entropy density $\eta=\eta(\rho,\vartheta)$.

      A general constraint identity \eqref{eq:OurAlgo:RunEx:ConstrID} yields the requirement
      \[
      \dfrac{\partial \eta}{\partial\vartheta} T_{12}=0,
      \]
      so the fluid must be isotropic, $T_{12}=0$, or the entropy must be temperature-independent. In the general case, one has
      \beq \label{eq:eg2:isotr}
      T_{12}=0,\qquad T_{11}=T_{22}=-\dfrac{\rho^2}{{\partial \eta}/{\partial\vartheta}}\left\{ \epsilon,\eta\right\}_{\left\{\rho,\vartheta\right\}} =: -P(\rho, \vartheta, \partial_i \vartheta),
      \eeq
      that is, the medium is \emph{isotropic}, and the diagonal components of the stress tensor can be associated with negative hydrostatic pressure $P$. The latter is given in terms of a Poisson bracket
      \beq\label{po:bra}
      \left\{ \epsilon,\eta\right\}_{\left\{\rho,\vartheta\right\}} = \dfrac{\partial \epsilon}{\partial\rho}\dfrac{\partial \eta}{\partial\vartheta}-\dfrac{\partial \epsilon}{\partial\vartheta}\dfrac{\partial \eta}{\partial\rho}.
      \eeq
      One can show that the entropy fluxes in the general case are defined by
      \beq \label{eq:eg2:Phis}
      \Phi_{1}(\rho, \vartheta, \partial_i \vartheta)=-A(\rho, \vartheta)\partial_2 \vartheta + B(\rho, \vartheta),\qquad
      \Phi_{2}(\rho, \vartheta, \partial_i \vartheta)=A(\rho, \vartheta)\partial_1 \vartheta + C(\rho, \vartheta),
      \eeq
      where $A$, $B$ and $C$ are arbitrary functions of their respective arguments.

      The remaining relationships define the energy fluxes $q_i=q_i(\rho, \vartheta, \partial_i \vartheta)$:
      \beq\label{eq:eg2:gen:Qs}
      \barr
      \dfrac{\partial  q_{1}}{\partial (\partial_1 \vartheta)} = \dfrac{\partial  q_{2}}{\partial (\partial_2 \vartheta)}=0,\qquad \dfrac{\partial  q_{1}}{\partial (\partial_2 \vartheta)} + \dfrac{\partial  q_{2}}{\partial (\partial_1 \vartheta)}=0,\qquad %\\[2ex]

      \dfrac{\partial  q_{i}}{\partial \rho} \dfrac{\partial  \eta}{\partial \vartheta}= \dfrac{\partial  \Phi_{i}}{\partial \rho} \dfrac{\partial  \epsilon}{\partial \vartheta},\quad i=1,2.
      \earr
      \eeq

\end{enumerate}

%\subsection{The comparison of the entropy principles}\label{sec:ep:comp}

%M-L: neither necessary nor sufficient...
%
%coincides with SolSet when stuff is linear, as shown in e.g. ***
%
%fails otherwise, as shown in e.g. ***
%
%leads to extra complications - Lambdas to be found... see e.g. LAST where eqs are complicated...

%\newpage
%-------------------------------------------------------------------------------------------------------------------------------------------------
\section{Computational example 1: One-dimensional gas dynamics}\label{sec:Maple:EG1}
%-------------------------------------------------------------------------------------------------------------------------------------------------
We now implement the solution set-based scheme, presented in Section \ref{sec:Entr:GenD}, in \verb|Maple|. For this, we resort to the symbolic software package \verb|GeM| for \verb|Maple|, which can be used for computations in conjunction with the symbolic manipulation of differential equations and their differential consequences, of related expressions, involving dependent and independent variables, partial derivatives and constitutive functions. This package incorporates routines for the computation of local conservation laws, as well as higher-order local symmetries of ordinary and partial differential equations and ODE/PDE systems \cite{GemReferenceOnline,cheviakov2007gem,cheviakov2010symbolic}. Furthermore, the package can be applied combined with standard \verb|Maple| routine \verb|rifsimp| for the reduction of overdetermined systems, for the analysis of conservation laws, or problems of symmetry. For details on the methods, routines and the representation of differential equations in \verb|GeM|, see the preceding work \cite{cheviakov2017symbolic} and \cite{cheviakov2007gem,cheviakov2010symbolic}.

The \verb|GeM| package is employed in the current work to facilitate the exploitation of the entropy principle via a newly proposed scheme, the solution set entropy principle. It automates computations that are exhausting on paper for more complex system and allows for easy manipulation of parameters and variable influences. Below, four examples are given to outline the implementation, starting with an introductory example of one-dimensional gas dynamics, followed by an example of a two-dimensional heat-conducting fluid, then a non-simple fluid for which differential consequences are of importance in the solution set approach and, finally, in Appendix \ref{sec:Maple:EG4}, a model of a granular flow.

%
%-------------------------------------------------------------------------------------------------------------------------------------------------
%---
%-------------------------------------------------------------------------------------------------------------------------------------------------
\subsection{The symbolic computation framework for the solution set approach}\label{sec:Maple:EG1:framew}
%{\color{red}***file:\verb|Eg1_gas_dynamics_1d_non-adiabatic_rhoW_SolSet_Dec03.mw|***}\\

In order to outline the symbolic computation algorithm, to explore the possibilities in the evaluation of the results and to perform a first comparison, we start with the rather elementary example of a (1+1)-dimensional system in the context of gas dynamics. The considered balance of mass, momentum and energy are
\begin{subequations}\label{eq:gas1D:gen0}
\begin{equation}\label{eq:gas1D:gen:dens0}
\Pi^{\rho}:~\partial_t \rho  +\partial_1 (\rho u)=0,
%\Pi^{v}: \rho D_t v - (\div T) =0,\qquad i=1\ldots, n,\\
\end{equation}
\begin{equation}\label{eq:gas1D:gen:momx0}
\Pi^{u}:~\rho (\partial_t u +u \partial_1 u ) - \partial_1 p =0, \\[2ex]
\end{equation}
\begin{equation}\label{eq:gas1D:gen:en0}
\Pi^{\epsilon}:~\rho (\partial_t \epsilon + u \partial_1 \epsilon) +  \partial_1 q_1  - p \partial_1 u =0,
\end{equation}
\end{subequations}
involving the physical fields given by the density $\rho=\rho(x,t)$, the velocity $u=u(x,t)$ in the $x-$direction, and  the internal energy $\epsilon=\epsilon(x,t)$. The model \eqref{eq:gas1D:gen0} contains the unknown constitutive functions
\beq\label{eq:qg1:const_f}
\psi = \left(p, q_1, \eta, \Phi_1 \right),
\eeq
i.e., the pressure $p$, the $x$-component of the energy flux $q_1$, the entropy density $\eta$, and the $x$-component of the entropy flux $\Phi_1$. The constitutive functions are assumed to depend on a simple constitutive class
\beq\label{eq:qg1:const_dep}
\phi_C= (\rho, \epsilon).
\eeq
One can equivalently use the temperature $\vartheta$ as a state variable instead of the internal energy $\epsilon$; this leads to additional equations. We note that instead of a more general stress tensor, the current inviscid gas dynamics model involves the scalar hydrostatic pressure $p$.

The one-dimensional entropy inequality is given by
\begin{equation}\label{eq:entropy_ineq:1Dgas0}
\Pi^{\eta}:~\rho (\partial_t \eta + u \partial_1 \eta ) +  \partial_1 \Phi_1 \geq 0.
\end{equation}

We now apply the solution set-based entropy principle to the model governed by Eqs. \eqref{eq:gas1D:gen0}-\eqref{eq:entropy_ineq:1Dgas0}, using the appropriate \verb|Maple| commands as shown below.

%%%%%%%%%%%%%%%%%%
\medskip\noindent \textbf{Step A. Initialize.} Clear the variables. Initialize the \verb|GeM| package.
\beq\label{eq:gem:initeq}
\barr
\verb|restart;|\\
\verb|read("d:/gem32_12.mpl");|\\
\earr
\eeq

%%%%%%%%%%%%%%%%%%
\medskip\noindent \textbf{Step B. Declare variables and constitutive functions.} The second step is to define the independent variables, the dependent variables (fields), and the constitutive functions, and declare them as the respective class of variables in \verb|Maple|/\verb|GeM|.
\[
\barr
\verb|ind:=t,x;|\\
\verb|dep:=R(ind), U(ind), E(ind);|\\
\verb|Constit_Dependence:=R(ind),E(ind);|\\
\verb|Constit_F:=P(Constit_Dependence), Q1(Constit_Dependence),|\\
\verb|           S(Constit_Dependence), Phi1(Constit_Dependence);|\\[2ex]
\verb|gem_decl_vars(indeps=[ind], deps=[dep], freefunc=[Constit_F]);|
\earr
\]

In \verb|Maple|, we apply the notation $\rho=\verb|R|$, $u=\verb|U|$, $\epsilon=\verb|E|$, $q_1=\verb|Q1|$, $\eta=\verb|S|$, $\Phi_1=\verb|Phi1|$, $p=\verb|P|$, for the fields and constitutive functions. It is our common convention to use small letters for independent variables, and capitals for dependent variables and constitutive functions.

%%%%%%%%%%%%%%%%%%
\medskip\noindent \textbf{Step C. Declare the model equations.} A symbolic operator for the material derivative $\sg{D}_t$, see Eq. \eqref{eq:TotD}, is introduced:
\beq\label{maple:Dt}
\verb|MaterialDer:=Q->diff(Q,t)+U(ind)*diff(Q,x);|\\
\eeq
The model PDEs \eqref{eq:gas1D:gen0} are now defined and, with the last command, declared in terms of \verb|Maple| symbols:
\beq\label{eq:sec31:defeqs}
\barr
\verb|Pi_Rho:=diff(R(ind),t) + diff(R(ind)*U(ind),x) = 0;|\\
\verb|Pi_U:=R(ind)*MaterialDer(U(ind)) + diff(P(Constit_Dependence),x) = 0;|\\
\verb|Pi_E:=R(ind)*MaterialDer(E(ind)) + diff(Q1(Constit_Dependence),x)|\\
\verb|      +P(Constit_Dependence)*diff(U(ind),x) = 0;|\\[2ex]
\verb|gem_decl_eqs( [Pi_Rho, Pi_U, Pi_E] );|
\earr
\eeq
This symbolic representation of the PDEs can be extracted using \verb|GeM| variables as follows:
\beq\label{eq:sec31:defeqs2}
\barr
\verb|Pi_R_Symb:=GEM_ALL_EQ_AN[1];|\\
\verb|Pi_U_Symb:=GEM_ALL_EQ_AN[2];|\\
\verb|Pi_E_Symb:=GEM_ALL_EQ_AN[3];|
\earr
\eeq
For example, the continuity equation $\Pi^{\rho}$ \eqref{eq:gas1D:gen:dens0} takes the symbolic form
\[
\verb|Pi_R_Symb = R*Ux+Rx*U+Rt|
\]

%%%%%%%%%%%%%%%%%%
\medskip\noindent \textbf{Step D. Leading derivatives; entropy inequality on solutions.} Choosing the time derivatives $\partial_t \rho$, $\partial_t u$ and $\partial_t \epsilon$ as the leading derivatives, the PDEs are written in the solved form, as follows:
\beq\label{eq:gasEG:leadDer}
\barr
\verb|Leading_ders:={Rt,Ut,Et};|\\
\verb|solved_DEs:=solve([Pi_R_Symb,Pi_U_Symb,Pi_E_Symb], Leading_ders);|
\earr
\eeq
Furthermore, the left-hand side of the entropy inequality \eqref{eq:entropy_ineq:1Dgas0} is defined and then converted into a symbolic \verb|Maple| expression with the following \verb|Maple|/\verb|GeM| commands.
\[
\barr
\verb|Pi_S:= R(ind)*MaterialDer(S(Constit_Dependence))|\\
\verb|      + diff(Phi1(Constit_Dependence),x);|\\
\verb|Pi_S_Symb:=gem_analyze(Pi_S);|
\earr
\]
In order to compute the \emph{entropy inequality on solutions}, the relations between derivatives arising from the system equations should be substituted into the entropy inequality \eqref{eq:entropy_ineq:1Dgas0}. Since the latter does not involve any derivatives of the leading derivatives $\partial_t \rho$, $\partial_t u$ and $\partial_t \epsilon$, no differential consequences are required; only the solved forms of the model equations \verb|solved_DEs| must be substituted into the entropy inequality:
\beq\label{eq:GesEG:entropyIneqSS}
\barr
\verb|Entropy_Inequality_Sol_Set:= simplify(subs(solved_DEs, Pi_S_Symb));|
\earr
\eeq
The entropy inequality on the solution set for the current example consequently becomes
\beq\label{eq:eg1:gasdyn:entr_ineq}
-\rho^2 \partial_1u \dfrac{\partial \eta}{\partial\rho}-\left(\partial_1\epsilon \dfrac{\partial q_1}{\partial\epsilon} + p \partial_1u + \partial_1\rho \dfrac{\partial q_1}{\partial\rho}\right)\dfrac{\partial \eta}{\partial\epsilon}   + \partial_1\epsilon \dfrac{\partial \Phi_1}{\partial\epsilon}+\partial_1\rho \dfrac{\partial \Phi_1}{\partial\rho}\geq 0.
\eeq

%%%%%%%%%%%%%%%%%%
\medskip\noindent \textbf{Step E. Constraint equations and residual entropy inequality.} All partial derivatives and the constitutive dependence of the problem in the symbolic form are obtained using the commands
\beq\label{eg1:step:E:1}
\barr
\verb|All_derivatives:={seq(GEM_ALL_ORDER_DERS[i][], i=1..GEM_MAX_ORDER)};|\\
\verb|Constit_Dependence_Symb:={map(x->gem_analyze(x), |\\
\qquad \qquad \verb| [Constit_Dependence])[]};|
\earr
\eeq
In particular, we have \verb|All_derivatives={Et, Ex, Rt, Rx, Ut, Ux}| and \verb|Constit_Dependence_Symb = {E, R}| in this example.

The \emph{free elements} of the current example include all independent and dependent variables and all derivatives present, with the exception of the leading derivatives and the constitutive dependence:
\beq\label{eg1:step:E:2}
\barr
\verb|Arbitrary_Elements:= {GEM_INDEP_VARS[], GEM_DEP_VARS[]}|\\
\verb|                     union All_derivatives|\\
\verb|                     minus Leading_ders|\\
\verb|                     minus Constit_Dependence_Symb;|
\earr
\eeq
The entropy inequality on the solution set \eqref{eq:eg1:gasdyn:entr_ineq} has to hold for all values of free elements, that is, for all solutions of the gas dynamics model  \eqref{eq:gas1D:gen0}.

In principle, as we will see in the examples below, the entropy inequality on the solution set  can have a denominator that must not vanish, whereas the numerator will produce the constraint equations and the residual entropy inequality. The numerator and the denominator are obtained using the commands
\beq\label{eg1:step:E:3}
\barr
\verb|Entropy_Inequality_numer:=expand(numer(Entropy_Inequality_Sol_Set));|\\
\verb|Entropy_Ineqality_denom:=denom(Entropy_Inequality_Sol_Set);|
\earr
\eeq
The preliminary constraint equations are obtained by setting the coefficients of the entropy inequality numerator to zero:
\beq\label{eg1:step:E:4}
\barr
\verb|coeffs_constraints:=coeffs(Entropy_Inequality_numer, Arbitrary_Elements);|
\earr
\eeq
The free coefficient (the part of the entropy inequality on the solution set that does not involve any free elements) yields the left-hand side of residual entropy inequality:
\beq\label{eg1:step:E:residual}
\barr
\verb|Residual_Entropy_Ineq:=simplify(eval(subs(|\\
\verb|    map(x->x=0,Arbitrary_Elements), Entropy_Inequality_numer)))|\\
\verb| /Entropy_Ineqality_denom;|
\earr
\eeq
It is identically zero in the current example. Finally, the actual constraint equations of the solution set method are obtained by the command
\beq\label{eg1:step:E:5}
\barr
\verb|Solution_Set_Constraints:=convert({coeffs_constraints} | \\
\verb|      minus {numer(Residual_Entropy_Ineq)}, list);|
\earr
\eeq
The constraint equations for the current example are given by three PDEs
\beq\label{eq:eg1:SolSet Meth} %Eg1_gas_dynamics_1d_non-adiabatic_rhoW_SolSet.mw
\barr
 \dfrac{\partial \Phi_1}{\partial \rho} - \dfrac{\partial \eta}{\partial \epsilon} \dfrac{\partial q_1 }{\partial \rho}=0,\qquad  \dfrac{\partial \Phi_1 }{\partial \epsilon} - \dfrac{\partial \eta}{\partial \epsilon} \dfrac{\partial q_1 }{\partial \epsilon}=0,\qquad
  p\dfrac{\partial \eta}{\partial \epsilon}+ \rho^2 \dfrac{\partial \eta}{\partial \rho}=0.
\earr
\eeq

%
%--------------------------------------------------------
\subsection{Results, comparisons, and discussion }

If the usual M\"{u}ller-Liu procedure (described in  Section \ref{sec:MuellerLiuAppr}) is applied to the same model, with Lagrange multipliers depending also on \eqref{eq:qg1:const_dep}, one obtains the following set of \emph{Liu identities}, which are deduced in more details in \cite{cheviakov2017symbolic,muller1998rational}:
\beq\label{eq:eg1:LiuMeth} %Eg1_gas_dynamics_1d_non-adiabatic_rhoW_Liu.mw
\barr
\rho \Lambda^{u}=0, \qquad \rho\Lambda^{\rho} +  \rho u \Lambda^{u} + p\Lambda^{\epsilon}=0, \\[2ex]

\rho \dfrac{\partial \eta}{\partial \rho} - \Lambda^{\rho}=0, \qquad \rho \left(\dfrac{\partial \eta}{\partial \epsilon} - \Lambda^{\epsilon}\right)=0,\\[2ex]

 \rho u\left( \dfrac{\partial \eta}{\partial \epsilon} - \Lambda^{\epsilon}\right)- \Lambda^{u}\dfrac{\partial p}{\partial \epsilon}
 - \Lambda^{\epsilon} \dfrac{\partial q_1}{\partial \epsilon} + \dfrac{\partial \Phi_1}{\partial \epsilon} =0,\\[2ex]

 v\left(\rho \dfrac{\partial \eta}{\partial \rho} - \Lambda^{\rho}\right) - \Lambda^{u}\dfrac{\partial p}{\partial \rho}
  - \Lambda^{\epsilon} \dfrac{\partial q_1}{\partial \rho} + \dfrac{\partial \Phi_1}{\partial \rho}  =0,
\earr
\eeq
for the constitutive functions $p, q_1, \eta, \Phi_1$ \eqref{eq:qg1:const_f} and the Lagrange multipliers $\Lambda^{\rho}$,  $\Lambda^{u}$, $\Lambda^{\epsilon}$ that correspond to the mass, momentum, and energy balance equations \eqref{eq:gas1D:gen0}. The PDEs \eqref{eq:eg1:LiuMeth} can be separated into two groups: one for the Lagrange multipliers, which yields
\beq\label{eq:eg1:LiuMeth:Lambdas}
\Lambda^{u} = 0,\qquad  \Lambda^{\rho}= \rho \dfrac{\partial \eta}{\partial \rho}, \qquad  \Lambda^{\epsilon}= \dfrac{\partial \eta}{\partial \epsilon},
\eeq
and another group for the physical constitutive functions, which gives the PDEs \eqref{eq:eg1:SolSet Meth}.

\medskip Comparing the results of the solution set-based method and the M\"{u}ller-Liu procedure, we observe that the constraint equations on the constitutive functions \eqref{eq:qg1:const_f} obtained by both methods indeed coincide. The solution set-based approach in this case offers a computational shortcut, since the calculations and the constraint equations \eqref{eq:eg1:SolSet Meth} involve no multipliers upfront.

\begin{remark}
In the current example, the constitutive dependence \eqref{eq:qg1:const_dep} is `simple', in the sense of not involving the leading derivatives \eqref{eq:gasEG:leadDer}. As a result, the extended entropy inequality \eqref{eq:GesEG:entropyIneqSS} computed on solutions of the given system does not require the substitution of differential consequences of the balance equations \eqref{eq:gas1D:gen0}, and the entropy inequality remains linear in the highest derivatives. This is the reason why
the results of the solution set method and the M\"{u}ller-Liu procedures are the same. This will not be the case for more complex constitutive dependencies, as it will be shown below.
\end{remark}

\begin{remark} It is of interest to consider the form of the entropy inequality for the well-known specific case of an \emph{ideal gas}. For a thermally and calorically perfect gas, the thermal energy density and the temperature are related by
\beq\label{eq:ideal:EC}
\epsilon=C_v\theta,
\eeq
where $C_v = {\tilde{R}}/{(\gamma-1)}$ is the specific heat at constant volume, $\tilde{R}=R/M$ is the specific gas constant, and $R$ and $M$ are the universal gas constant and the molar mass of the gas. Moreover, for an ideal gas, one has $p=\rho \hat{R} \,\theta$, $\eta = C_v \ln (\rho^{1-\gamma}\,\theta)$, so the specific forms of the constitutive functions $p, \eta$ in \eqref{eq:qg1:const_f} are given by
\beq\label{eq:ideal:PS}
p(\rho, \epsilon)=(\gamma-1) \rho \epsilon,\qquad \eta(\rho, \epsilon) = C_v \ln \left(\dfrac{\epsilon}{C_v\, \rho^{\gamma-1}}\right),
\eeq
%and the Lagrange multipliers of the M\"{u}ller-Liu procedure obtained using \eqref{eq:eg1:LiuMeth:Lambdas} are
%\[
%\Lambda^{\rho}(\rho, \epsilon) = - \hat{R}, \qquad \Lambda^{v}(\rho, \epsilon)=0,\qquad \Lambda^{\epsilon}(\rho, \epsilon)=\dfrac{C_v}{\epsilon}=\dfrac{1}{\theta}.
%\]
and as per \eqref{eq:eg1:SolSet Meth}, the energy and entropy fluxes $q_1$ and $\epsilon$ satisfy
\beq\label{eq:ideal:QEps:gen}
\dfrac{\partial \Phi_1}{\partial \rho} = \dfrac{C_v}{\epsilon} \dfrac{\partial q_1 }{\partial \rho}=0,\qquad  \dfrac{\partial \Phi_1 }{\partial \epsilon} - \dfrac{C_v}{\epsilon} \dfrac{\partial q_1 }{\partial \epsilon}=0.
\eeq
The entropy inequality $\Pi^{\eta}$ \eqref{eq:entropy_ineq:1Dgas0} in this case becomes
\beq\label{eq:ideal:entr}
\left((\gamma-1)\dfrac{\partial \Phi_1}{\partial \rho} - \dfrac{\tilde{R}}{\epsilon}\dfrac{\partial q_1}{\partial \rho} \right)\partial_1 \rho + \left((\gamma-1)\dfrac{\partial \Phi_1}{\partial \epsilon} - \dfrac{\tilde{R}}{\epsilon}\dfrac{\partial q_1}{\partial \epsilon} \right)\partial_1 \epsilon\geq 0.
\eeq
A particular solution of \eqref{eq:ideal:QEps:gen} is given by vanishing entropy and energy fluxes, $\Phi_1=q_1=0$, whereupon the entropy inequality \eqref{eq:ideal:entr} vanishes identically. Since the entropy flux is  zero, this corresponds to  locally adiabatic gas motion, with the vanishing entropy production rate $\frac{d}{dt} \int_{I(t)} \rho \eta \, dx =0$ in any material interval $I(t)$ (cf. \eqref{eq:glob:Q2}).
%
%d/dt ???dV (dx??)=0 for every V... see intro...******************
%
%
%
% +\oint_{\partial V(t)} \Phi_i dS_i - \int_{V(t)} \rho s \, dV.
%
%
%For a locally adiabatic gas motion, i.e., no interaction with the environment,

\end{remark}

%\medskip\noindent \textbf{Step F. On the results}

%
%--------------------------------------------------------
\subsection{Case splitting for the constraint equations}\label{sec:eg1:splitt}

The entropy inequality as in \eqref{eq:entropy_ineq:Liufluid} or \eqref{eq:entropy_ineq:1Dgas0}, can be stated for a class of materials when there are grounds to believe that the entropy can indeed be defined for each model in a given class as a function of the state, that is, in terms of the state variables. For every model in the class, the entropy density $\eta$ and the entropy fluxes $\Phi_i$ are given by some specific expressions, which may be \emph{a priori} unknown to a researcher.

The associated constraint equations, like the PDEs \eqref{eq:OurAlgo:RunEx:ConstrID} or \eqref{eq:eg1:SolSet Meth}, are nonlinear and underdetermined; they may provide different restrictions on the form of the remaining constitutive functions for different forms of the entropy-related constitutive functions $\eta$, $\Phi_i$. It may therefore be useful to employ the available mathematical and computational tools that enable to split and analyze the constraint equations with respect to possible forms of $\eta$ and $\Phi_i$. In particular, the Gr\"obner basis-based technique implemented in the \verb|Maple| \verb|rifsimp| routine can be used for such case splitting.

We illustrate the case splitting idea with the aid of the current example, for the given class of gas dynamics models \eqref{eq:gas1D:gen0} and the corresponding constraint equations \eqref{eq:eg1:SolSet Meth}. The aim is to find out what types of constraints on the forms of the constitutive functions $p(\rho, \epsilon), q_1(\rho, \epsilon)$ one may have, based on different functional forms of $\eta(\rho, \epsilon)$ and $\Phi_1(\rho, \epsilon)$. The case splitting is implemented in the \verb|rifsimp|/\verb|casesplit| command, and the following code is used:
\beq\label{eq:gas:rss}
\barr
\verb|Split_WRT_entropy_quantities:=DEtools[rifsimp](|\\
\verb|    coeffs_constraints, [P(R, E), Q1(R, E)],|\\
\verb|    mindim=1, casesplit);|
\earr
\eeq
In particular, the arguments $p(\rho, \epsilon), q_1(\rho, \epsilon)$ provided to \verb|rifsimp| are the unknown functions, and other functions present in the constraint equations \verb|coeffs_constraints| (in this case, $\eta(\rho, \epsilon)$ and $\Phi_1(\rho, \epsilon)$) are treated by \verb|rifsimp| as classifying functions, whose different forms may yield different solution classes (cases) for the constitutive functions $p(\rho, \epsilon), q_1(\rho, \epsilon)$. The output variable \verb|Split_WRT_entropy_quantities| is a Maple table, containing four different cases. The resulting \emph{case tree} can be plotted using the command:
\beq\label{eq:gas:caseplot}
\verb|DEtools[caseplot](Split_WRT_entropy_quantities, pivots);|
\eeq

\begin{figure}[!h]
  \begin{center}
  \includegraphics[width=0.4\textwidth]{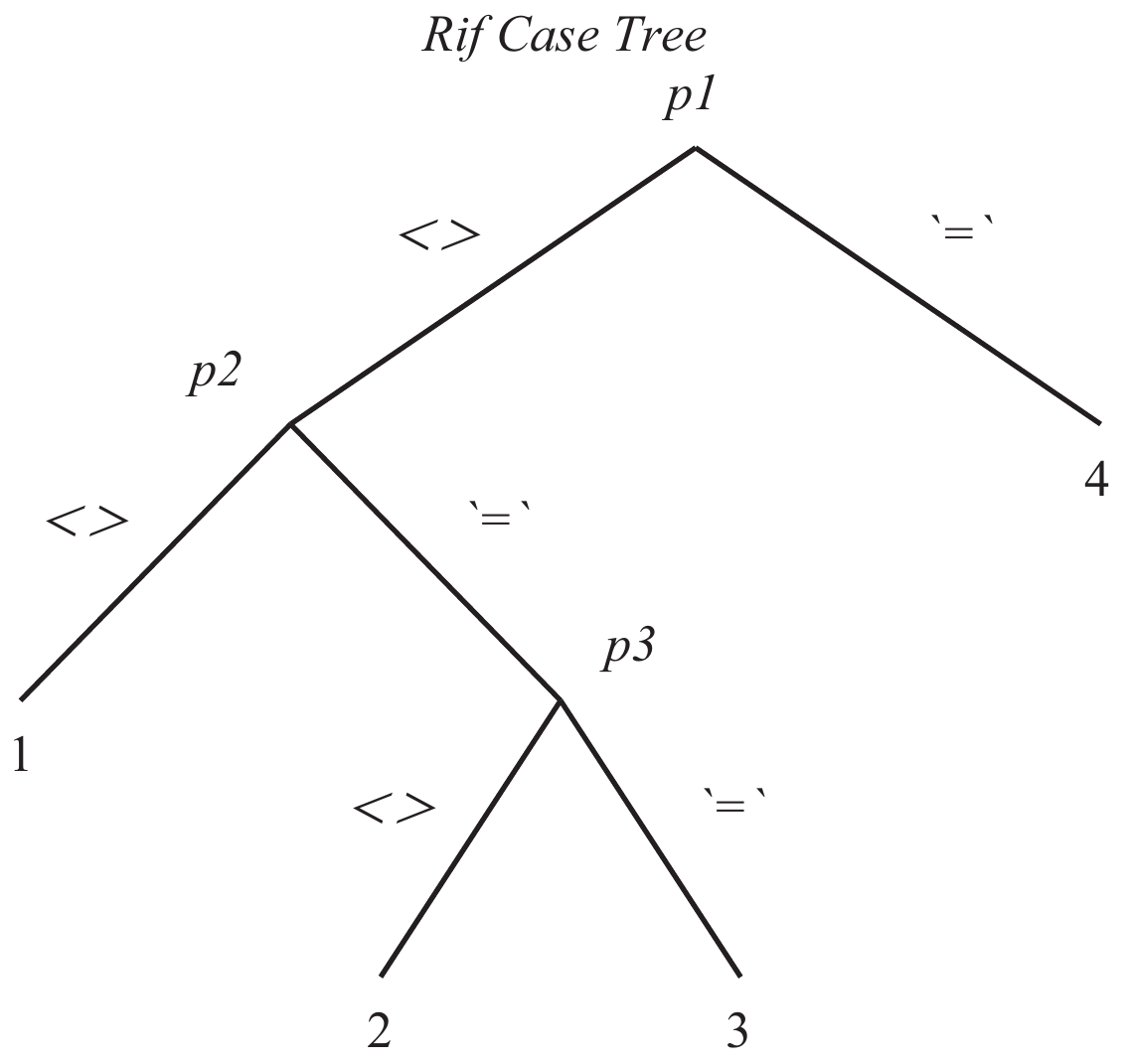}\\
  \end{center}
  \caption{Case tree for the classification of restrictions on the constitutive functions $p(\rho, \epsilon)$, $q_1(\rho, \epsilon)$ with respect to the forms of the entropy-related constitutive functions $\eta(\rho, \epsilon)$, $\Phi_1(\rho, \epsilon)$, in the framework of the solution set entropy principle applied to the one-dimensional gas dynamics model \eqref{eq:gas1D:gen0}. The constitutive functions and the constitutive class are given by \eqref{eq:qg1:const_f}, \eqref{eq:qg1:const_dep}.}\label{fig:1Dgasnonadia:casetree}
\end{figure}

The tree contains three pivots and four case branches (Figure \ref{fig:1Dgasnonadia:casetree}); the pivots, which are either equal to zero or not, are given by
\beq\label{eq:gasdyn1d:pivots}
\verb|p1|=\frac{\partial \eta}{\partial \epsilon},\quad \verb|p2|=\frac{\partial \Phi_1}{\partial \epsilon},\quad \verb|p3|=\frac{\partial \Phi_1}{\partial \rho}.
\eeq
We now analyze the four cases that arise.

\medskip\noindent\textbf{Case 1.} This is the most general case given by \eqref{eq:eg1:SolSet Meth}, with no additional restrictions on the forms of $\eta(\rho, \epsilon)$ and $\Phi_1(\rho, \epsilon)$, except for a condition
\beq\label{eq:gas:SPh:case1}
\dfrac{\partial \Phi_1 }{\partial \epsilon}\dfrac{\partial^2 \eta}{\partial\rho\partial \epsilon} =
\dfrac{\partial \Phi_1 }{\partial \rho}\dfrac{\partial^2 \eta}{\partial \epsilon^2},
\eeq
following from the first two PDEs of \eqref{eq:eg1:SolSet Meth}. The condition \eqref{eq:gas:SPh:case1} is a zero Poisson bracket $\{\Phi_1, {\partial \eta}/{\partial \epsilon}\}_{\left\{\rho,\epsilon\right\}}=0$, hence these two functions are functionally dependent: there exists a function $F(A,B)$ so that
\beq\label{eq:gas:SPh:case1:b}
F(\Phi_1,{\partial \eta}/{\partial \epsilon})=0,
\eeq
or, where it can be solved for $\Phi_1$, one has
\beq\label{eq:gas:SPh:case1:c}
\Phi_1(\rho, \epsilon)=G({\partial \eta}/{\partial \epsilon}),
\eeq
for an arbitrary sufficiently smooth function $G=G(z)$ of one variable.

For every pair of functions $\eta(\rho, \epsilon)$, $\Phi_1(\rho, \epsilon)$ satisfying \eqref{eq:gas:SPh:case1:b} or  \eqref{eq:gas:SPh:case1:c}, the other two constitutive functions, the pressure $p(\rho, \epsilon)$ and the energy flux $q_1(\rho, \epsilon)$, are determined from
\beq\label{eq:gas:QP:case1}
p(\rho, \epsilon)=- \rho^2 \dfrac{{\partial \eta}/{\partial \rho}}{{\partial \eta}/{\partial \epsilon}},\qquad
\dfrac{\partial q_1 }{\partial \rho} = \dfrac{{\partial \Phi_1}/{\partial \rho}}{{\partial \eta}/{\partial \epsilon}}, \qquad
\dfrac{\partial q_1 }{\partial \epsilon} = \dfrac{{\partial \Phi_1}/{\partial \epsilon}}{{\partial \eta}/{\partial \epsilon}}.
\eeq
In particular, the second and third equations of \eqref{eq:gas:QP:case1} yield another vanishing Poisson bracket $\{q_1, \Phi_1\}_{\left\{\rho,\epsilon\right\}}=0$, and consequently, the pressure and the energy flux in Case 1 are given by
\beq\label{eq:gas:QP:case1:b}
p(\rho, \epsilon)=- \rho^2 \dfrac{{\partial \eta}/{\partial \rho}}{{\partial \eta}/{\partial \epsilon}},\qquad
q_1(\rho, \epsilon) = H({\partial \eta}/{\partial \epsilon}),
\eeq
where $H=H(z)$ satisfies $H'(z)=G'(z)/z$.

\medskip\noindent\textbf{Case 2.} Here, with the pivot \verb|p2|$=0$, one has the reduced dependence $\Phi_1=\Phi_1(\rho)$ only. The constraint equations \eqref{eq:eg1:SolSet Meth} in this case yield a simple solution
\beq\label{eq:eg1:Casespl:Case2}
\barr
\Phi_1(\rho, \epsilon) = f_1(\rho),\qquad \eta(\rho, \epsilon) = f_2(\rho)\epsilon +f_3(\rho),\\[3ex]
q_1(\rho, \epsilon) = \displaystyle\int \dfrac{f_1'(\rho)}{f_2(\rho)}\, d\rho,\qquad  p(\rho, \epsilon)=-\rho^2\dfrac{f_2'(\rho)\epsilon +f_3'(\rho)}{f_2(\rho)},
\earr
\eeq
where $f_1, f_2, f_3$ are arbitrary functions of $\rho$. It is evident that Case 2 is just a special case of Case 1, considered for a particular form of $\Phi_1$, and yields no new information.

\medskip\noindent\textbf{Case 3.} In this case, both pivots \verb|p2|=\verb|p3|$=0$, and one has $\Phi_1(\rho, \epsilon)=c_1=\const$. The remaining  constraint equations \eqref{eq:eg1:SolSet Meth} yield
\beq\label{eq:eg1:Casespl:Case3}
\barr
\eta = \eta(\rho, \epsilon),\qquad q_1(\rho, \epsilon) =c_2=\const, \qquad p(\rho, \epsilon)=-\rho^2\dfrac{{\partial \eta}/{\partial \rho}}{{\partial \eta}/{\partial \epsilon}}.
\earr
\eeq
It follows that, if the entropy flux density is constant, the energy flux must also be constant; moreover, in this case, there are no restrictions on the form of the entropy density. We note that the adiabatic gas motion, where $q_1=\Phi_1=0$, belongs to this case. Case 3 is not included in Case 1, and is new.

\medskip\noindent\textbf{Case 4.} Finally, in Case 4, the pivot \verb|p1| in \eqref{eq:gasdyn1d:pivots} is zero. The  constraint equations \eqref{eq:eg1:SolSet Meth} then yield
\beq\label{eq:eg1:Casespl:Case4}
\barr
\eta(\rho, \epsilon) = c_1=\const, \qquad \Phi_1(\rho, \epsilon) =c_2=\const.
\earr
\eeq
This is a degenerate special case, with the entropy density and flux being independent of the state of the model given by $\epsilon$ and $\rho$. Due to this, one has no restrictions on the forms of the gas pressure and the energy flux:
\beq\label{eq:eg1:Casespl:Case4:b}
\barr
p = p(\rho, \epsilon),\qquad q_1=q_1(\rho, \epsilon)
\earr
\eeq
are arbitrary. Formally, this case is also a new case, not contained in the previous ones.

\begin{remark}
In this classification example, we observe the common situation that more restrictive cases on the classifying functions $\eta(\rho, \epsilon)$, $\Phi_1(\rho, \epsilon)$ yield less restrictive requirements on the rest of the constitutive functions.
\end{remark}

\begin{remark}
We note that other classifications, that is, classifications with respect to other constitutive functions present in the problem, can be done. Such classifications would be equivalent to the one presented above, yet their forms of appearance can lead to further insights into possible relationships between the constitutive functions of this model.
\end{remark}

%-------------------------------------------------------------------------------------------------------------------------------------------------
%-------------------------------------------------------------------------------------------------------------------------------------------------
\section{Computational example 2: Two-dimensional heat-conducting fluid}\label{sec:eg2}
%
%{\color{red}***file:\verb|Eg2_SimpFluid_SolSetApproach_01, Eg2_SimpFluid_LiuApproach_01|***}\\

As a more extensive example, we now use the solution set approach to derive restrictions on constitutive functions for a heat-conducting fluid flow, as described in Section \ref{sec:Entr:GenD}, in the case of two spatial dimensions. The continuity, momentum, and energy equations, in the absence of external supply terms, are given by
\begin{subequations}\label{eq:fluid:2d:heatcond}
\begin{equation}\label{eq:fluid:2d:dens}
\Pi^{\rho}:~\partial_t \rho  +\partial_1 (\rho u)+ \partial_2(\rho v)=0,
%\Pi^{v}: \rho D_t v - (\div T) =0,\qquad i=1\ldots, n,\\
\end{equation}
\begin{equation}\label{eq:fluid:2d:momx}
\Pi^{u}:~\rho (\partial_t u +u \partial_1 u + v \partial_2 u) - \partial_1 T_{11} - \partial_2 T_{12} =0, \\[2ex]
\end{equation}
\begin{equation}\label{eq:fluid:2d:momy}
\Pi^{v}:~\rho (\partial_t v +u \partial_1 v + v \partial_2 v) - \partial_1 T_{12} - \partial_2 T_{22} =0, \\[2ex]
\end{equation}
\begin{equation}\label{eq:fluid:2d:en}
\Pi^{\epsilon}:~\rho (\partial_t \epsilon + u \partial_1 \epsilon+ v \partial_2 \epsilon ) +  \partial_1 q_1 +  \partial_2 q_2 - T_{11} \partial_1 u - T_{12} \partial_2 u- T_{12} \partial_1 v- T_{22} \partial_2 v=0,
\end{equation}
\end{subequations}
where $\rho(t,x,y)$ is the density, and $u(t,x,y)$, $v(t,x,y)$ are the scalar spatial velocity components, respectively. Instead of the hydrostatic pressure, three independent components of the stress tensor emerge, $T_{11}(t,x,y), T_{12}(t,x,y), T_{22}(t,x,y)$. The heat flux vector components are given by $q_1(t,x,y)$ and $q_2(t,x,y)$. The internal energy density is denoted by $\epsilon(t,x,y)$, and the temperature by $\vartheta(t,x,y)$.

The system of equations \eqref{eq:fluid:2d:heatcond} involves nine unknowns and four equations, and therefore is underdetermined. The closure for \eqref{eq:fluid:2d:heatcond} involves the constitutive relations and/or additional PDEs. We now apply the solution-set based entropy principle to derive some constitutive relations.

The entropy inequality in two dimensions is given by
\begin{equation}\label{eq:entropy_ineq:1Dgas}
\Pi^{\eta}:~\rho (\partial_t \eta + u \partial_1 \eta +v \partial_2 \eta) +  \partial_1 \Phi_1 +  \partial_2 \Phi_2 \geq 0,
\end{equation}
with the entropy density $\eta(t,x,y)$ and the entropy flux components $\Phi_1(t,x,y),\Phi_2(t,x,y)$. We let $\rho$, $u$, $v$, $\vartheta$ be the dependent variables, and assume the constitutive dependencies \eqref{eq:simplefluid:constit:our:algo} for the nine constitutive functions \eqref{eq:simplefluid:constitF:our:algo}.

%for $\psi = \left(T_{11},T_{12},T_{2}, \epsilon,q_1, q_2,\eta, \Phi_1,\Phi_2 \right)$ as in equation \eqref{eq:simplefluid:constit:our:algo}, with
%\begin{equation}
%\begin{split}
%\epsilon=\epsilon(\rho,\vartheta), \quad \eta=\eta(\rho,\vartheta),\\
%T_{ij}=T_{ij}(\rho, \vartheta, \partial_i \vartheta), \quad
% q_i=q_i(\rho, \vartheta, \partial_i \vartheta), \quad \Phi_{i}=\Phi_{i}(\rho, \vartheta, \partial_i \vartheta), \quad i,j=1,2.
%\end{split}
%\label{eq:gas1D:gen:constitutive}
%\end{equation}

\subsection{Symbolic computations for the solution set approach}\label{sec:eg2:comp}

The symbolic computation follows the same pattern as the previous computational example.

%%%%%%%%%%%%%%%%%%
\medskip\noindent \textbf{Step A. Initialize.} Initialization via commands \eqref{eq:gem:initeq}.

%%%%%%%%%%%%%%%%%%
\medskip\noindent \textbf{Step B. Declare variables and constitutive functions.} We now declare the dependent and independent variables of the problem, and the two different dependencies, as per the formulas \eqref{eq:simplefluid:constit:our:algo}, for the two sets of constitutive functions.
\beq\label{eg2:symb:setup}
\barr
\verb|ind:=t,x,y;   dep:=R(ind), U(ind), V(ind), W(ind);|\\
\verb|Constit_Dependence:=R(ind),W(ind),diff(W(ind),x),diff(W(ind),y);|\\
\verb|Constit_Dependence2:=R(ind),W(ind);|\\
\verb|Constit_F:=T11(Constit_Dependence), T12(Constit_Dependence),|\\
\verb|           T22(Constit_Dependence), Q1(Constit_Dependence), |\\
\verb|           Q2(Constit_Dependence), Phi1(Constit_Dependence),|\\
\verb|           Phi2(Constit_Dependence), |\\
\verb|           E(Constit_Dependence2), S(Constit_Dependence2);|\\
\earr
\eeq
\[
\barr
\verb|gem_decl_vars(indeps=[ind], deps=[dep], freefunc=[Constit_F]);|
\earr
\]
In \verb|Maple| notation, in addition to previously discussed quantities, we now have the second velocity component $v=\verb|V|$, the temperature $\vartheta=\verb|W|$, second flux components $q_2=\verb|Q2|$, $\Phi_2=\verb|Phi2|$, and the stress tensor components $T_{11}=\verb|T11|$, $T_{12}=\verb|T12|$, $T_{22}=\verb|T22|$ as constitutive variables. In order to implement the reduced dependencies of $\eta$ and $\epsilon$, a second constitutive class \verb|Constit_Dependence2| has been used.

%%%%%%%%%%%%%%%%%%
\medskip\noindent \textbf{Step C. Declare the model equations.} Using the two-dimensional symbolic operator for the material derivatives,
\beq\label{maple:Dt:2D}
\verb|MaterialDer:=Q->diff(Q,t)+U(ind)*diff(Q,x)+V(ind)*diff(Q,y);|\\
\eeq
the model PDEs \eqref{eq:fluid:2d:heatcond} can be defined:
\[
\barr
\verb|Pi_Rho:=diff(R(ind),t) + diff(R(ind)*U(ind),x)|\\
\verb|        + diff(R(ind)*V(ind),y) = 0;|\\
\verb|Pi_U:=R(ind)*MaterialDer(U(ind)) - diff(T11(Constit_Dependence),x)|\\
\verb|      - diff(T12(Constit_Dependence),y)=0; |\\
\verb|Pi_V:=R(ind)*MaterialDer(V(ind)) - diff(T12(Constit_Dependence),x)|\\
\verb|      - diff(T22(Constit_Dependence),y)=0; |\\
\verb|Pi_E:=R(ind)*MaterialDer(E(Constit_Dependence2))|\\
\verb|      + diff(Q1(Constit_Dependence),x)|\\
\verb|      + diff(Q2(Constit_Dependence),y) |\\
\verb|      - T11(Constit_Dependence)*diff(U(ind),x)|\\
\verb|      - T12(Constit_Dependence)*diff(U(ind),y)|\\
\verb|      - T12(Constit_Dependence)*diff(V(ind),x)|\\
\verb|      - T22(Constit_Dependence)*diff(V(ind),y) = 0;|\\
\earr
\]
The \verb|GeM| symbolic representation is obtained using the command
\[
\barr
\verb|gem_decl_eqs( [Pi_Rho, Pi_U, Pi_V, Pi_E] );|\\
\earr
\]
and extracted into four scalar variables
\beq\label{eq:get:Pi:1234}
\barr
\verb|Pi_R_Symb:=GEM_ALL_EQ_AN[1];|\\
\verb|Pi_U_Symb:=GEM_ALL_EQ_AN[2];|\\
\verb|Pi_V_Symb:=GEM_ALL_EQ_AN[3];|\\
\verb|Pi_E_Symb:=GEM_ALL_EQ_AN[4];|
\earr
\eeq
which now contain the symbolic forms of the model equations \eqref{eq:fluid:2d:heatcond}.

%%%%%%%%%%%%%%%%%%
\medskip\noindent \textbf{Step D. Leading derivatives; entropy inequality on solutions.} For the balance equations \eqref{eq:fluid:2d:heatcond}, it is natural to choose time derivatives  $\partial_t \rho$, $\partial_t u$, $\partial_t v$ and $\partial_t \epsilon$ as leading derivatives, since these are the highest-order derivatives only referring to a single independent variable. In order to describe the solution set of the two-dimensional heat-conducting fluid model, its governing equations \eqref{eq:fluid:2d:heatcond} are solved for the leading derivatives:
\[
\barr
\verb|Leading_ders:={Rt,Ut,Vt,Wt};|\\
\verb|solved_DEs:=solve([Pi_R_Symb,Pi_U_Symb,Pi_V_Symb,Pi_E_Symb],|\\
\verb|                  Leading_ders);|
\earr
\]
The left-hand side of the entropy inequality \eqref{eq:entropy_ineq:1Dgas} is now defined and converted into a symbolic  expression:
\beq\label{eq:eg2:PiS}
\barr
\verb|Pi_S:=R(ind)*MaterialDer(S(Constit_Dependence2))|\\
\verb|      + diff(Phi1(Constit_Dependence),x)|\\
\verb|      + diff(Phi2(Constit_Dependence),y);|\\
\verb|Pi_S_Symb:=gem_analyze(Pi_S);|
\earr
\eeq
The entropy inequality on the solution set is consequently computed, using the substitution \eqref{eq:GesEG:entropyIneqSS}. It is straightforward to verify that the result involves neither the leading derivatives nor their differential consequences.

%%%%%%%%%%%%%%%%%%
\medskip\noindent \textbf{Step E. Constraint equations and residual entropy inequality.} Next, the constraint equations are computed. We note that, even though we have two types of constitutive dependencies in this example, the one given by \verb|Constit_Dependence| includes the second one given by \verb|Constit_Dependence2|, therefore commands \eqref{eg1:step:E:1}--\eqref{eg1:step:E:5} work for this example without change.

The resulting constraint equations for this example are given by \eqref{eq:OurAlgo:RunEx:ConstrID}, and the residual entropy inequality by \eqref{eq:Liu:RunEx:ResIneq}.

\medskip\noindent \textbf{Case splitting.}
Commands similar to \eqref{eq:gas:rss} and \eqref{eq:gas:caseplot} can now be used to simplify the constraint equations, and split and plot cases, if any occur.
\beq\label{eq:eg2:casespl:plot}
\barr
\verb|Split_WRT_entropy_quantities:=DEtools[rifsimp](|\\
\verb|    [ Solution_Set_Constraints[], |\\
\verb|      Entropy_Ineqality_denom<>0,|\\
\verb|      diff(S(R, W),W)<>0, diff(E(R, W),W)<>0 ],|\\
\verb|    [ T11(R, W, Wx, Wy), T12(R, W, Wx, Wy), T22(R, W, Wx, Wy),|\\
\verb|      Q1(R, W, Wx, Wy), Q2(R, W, Wx, Wy),|\\
\verb|      Phi1(R, W, Wx, Wy), Phi2(R, W, Wx, Wy) ],|\\
\verb|    mindim=1, casesplit);|\\[2ex]
\verb|DEtools[caseplot](Split_WRT_entropy_quantities, pivots);|
\earr
\eeq
We chose the \emph{classifying functions} to be $\eta(\rho, \vartheta)$ and $\epsilon(\rho, \vartheta)$, hence they have been excluded from the set of unknown functions in the last vector-type argument of the \verb|rifsimp| command.

\begin{figure}[!h]
  \begin{center}
  \includegraphics[width=0.3\textwidth]{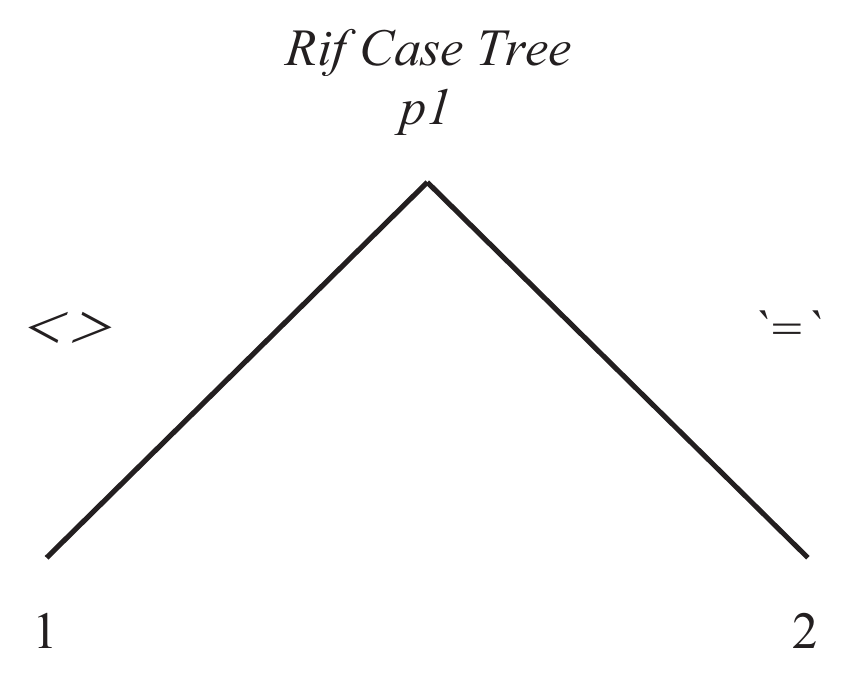}\\
  \end{center}
  \caption{The case tree for the classification of restrictions on the constitutive functions \eqref{eq:simplefluid:constitF:our:algo}, \eqref{eq:simplefluid:constit:our:algo} with respect to the classifying functions $\eta(\rho, \vartheta)$ and $\epsilon(\rho, \vartheta)$, for the two-dimensional heat-conducting fluid model \eqref{eq:fluid:2d:heatcond} (Section \ref{sec:eg2}). The pivot \emph{p1} is given by \eqref{eq:eg2:norm:pivot}.}\label{fig:eg2Gen}
\end{figure}

We observe that two cases arise (Figure \ref{fig:eg2Gen}), with the pivot given by
\beq\label{eq:eg2:norm:pivot}
\verb|p1|=\dfrac{\partial \epsilon}{\partial \vartheta} \dfrac{\partial^2 \eta}{\partial \rho\partial \vartheta} - \dfrac{\partial \eta}{\partial \vartheta} \dfrac{\partial^2 \epsilon}{\partial \rho\partial \vartheta}\, .
\eeq

\medskip\noindent\textbf{Case 1.} In the general Case 1, $\verb|p1|\ne 0$, the forms of the classifying functions $\epsilon=\epsilon(\rho,\vartheta)$ and $\eta=\eta(\rho,\vartheta)$ are arbitrary; the remaining constitutive function classes that correspond to nonnegative entropy production are defined by the constraint equations \eqref{eq:eg2:isotr}--\eqref{eq:eg2:gen:Qs}, and the residual entropy inequality is given by \eqref{eq:Liu:RunEx:ResIneq}.

In both Case 1 and Case 2, the isotropy condition $T_{12}=0$ holds, and the diagonal components of the stress tensor are associated with the hydrostatic pressure (see \eqref{eq:eg2:isotr}).

\medskip\noindent\textbf{Case 2.} This special case takes place when the pivot $\verb|p1|= 0$, which leads, for example, to the energy density being arbitrary, $\epsilon=\epsilon(\rho,\vartheta)$, and the entropy density given by
\beq\label{eq:eg2:case2}
\eta(\rho,\vartheta) = \int \left(\dfrac{\partial \epsilon}{\partial \vartheta} M(\vartheta)\right)\, d\vartheta + N(\rho),
\eeq
where $M$ and $N$ are arbitrary functions. As before, the reduced freedom in one of the classifying functions (the entropy density) leads to the increased freedom in other constitutive functions, such as the energy and entropy flux functions. The latter have to satisfy the PDEs
\beq \label{eq:eg2:PhiQ:cas2} % file Eg2_SimpFluid_SolSetApproach_01corr.mw
\barr
\dfrac{\partial \eta}{\partial \vartheta} \, \dfrac{\partial^2  q_{1}}{\partial (\partial_1 \vartheta)^2}= \dfrac{\partial \epsilon}{\partial \vartheta}\, \dfrac{\partial^2  \Phi_{1}}{\partial (\partial_1 \vartheta)^2}, \\[2ex]

\dfrac{\partial \eta}{\partial \vartheta}\left(\dfrac{\partial q_{2}}{\partial (\partial_1 \vartheta)} + \dfrac{\partial q_{1}}{\partial (\partial_2 \vartheta)} \right)=\dfrac{\partial \epsilon}{\partial \vartheta}\left(\dfrac{\partial \Phi_{2}}{\partial (\partial_1 \vartheta)} + \dfrac{\partial \Phi_{1}}{\partial (\partial_2 \vartheta)} \right), \\[2ex]

\dfrac{\partial  \eta}{\partial \vartheta}\, \dfrac{\partial  q_{i}}{\partial \rho} = \dfrac{\partial  \epsilon}{\partial \vartheta}\, \dfrac{\partial  \Phi_{i}}{\partial \rho} , \quad i=1,2,\\[2ex]

\dfrac{\partial  \eta}{\partial \vartheta}\,\dfrac{\partial  q_{i}}{\partial (\partial_i \vartheta)} = \dfrac{\partial  \epsilon}{\partial \vartheta}\, \dfrac{\partial  \Phi_{i}}{\partial ( \partial_i \vartheta)} , \quad i=1,2.\\[2ex]
\earr
\eeq
which are less restrictive than the constraint equations \eqref{eq:eg2:Phis}, \eqref{eq:eg2:gen:Qs} of the general Case 1. For example, here the energy and entropy flux functions $q_{i}, \Phi_{i}$ may depend on the corresponding gradient components of $\partial_i \vartheta$.

\begin{remark}
Case 1 and Case 2 above yield two different sets of constraints on the constitutive functions. Families of constitutive functions satisfying the restrictions of either Case 1 (constraint equations \eqref{eq:OurAlgo:RunEx:ConstrID}) or Case 2 (constraint equations \eqref{eq:eg2:PhiQ:cas2}) are consistent with the entropy principle, in the sense that they lead to a nonnegative entropy production inequality \eqref{eq:gen:EntrOnSol} on all solutions of the physical model \eqref{eq:fluid:2d:heatcond}, as long as that the residual entropy inequality \eqref{eq:Liu:RunEx:ResIneq} is satisfied.
\end{remark}

%
%--------------------------------------------------------------------------------------
\subsection{Results, comparisons, and discussion }\label{sec:eg2:compar}

For the physical model \eqref{eq:fluid:2d:heatcond} involving constitutive functions \eqref{eq:simplefluid:constit:our:algo}, the general constraint equations \eqref{eq:OurAlgo:RunEx:ConstrID} that follow from the solution set-based approach coincide with the constraint equations (Liu identities), computed in Ref.~\cite{cheviakov2017symbolic} using the classical M{\"u}ller-Liu procedure \cite{liu_method_1972}. The reason for it is that the constitutive dependencies \eqref{eq:simplefluid:constit:our:algo} are again  `simple', not involving the leading derivatives, and hence the entropy inequality computed on solutions remains linear in the highest derivatives. This is not the case for more complex constitutive dependencies, such as the one considered in Section \ref{sec:Maple:EG3} below, where the nonlinear terms will be present in the entropy inequality, and will not be accounted for by the M{\"u}ller-Liu approach.

\medskip
In the classical M{\"u}ller-Liu procedure, coefficients are set to zero only at higher-order derivatives which appear linearly in the extended entropy inequality \eqref{eq:ML:extEntr}. For the current example, the derivatives $\partial_i \vartheta$ are included in the constitutive dependence \eqref{eq:simplefluid:constit:our:algo}, and therefore are not used when setting to zero coefficients in the M{\"u}ller-Liu procedure; consequently, the residual entropy inequality \eqref{eq:Liu:RunEx:ResIneq} arises. By contrast, the solution set approach allows for the consideration of a nonlinear dependence (on higher-order derivatives) of the entropy inequality on solutions \eqref{eq:gen:EntrOnSol}. Indeed, here one could effectively include the residual entropy inequality \eqref{eq:Liu:RunEx:ResIneq} in the set of constraint identities given by \eqref{eq:OurAlgo:RunEx:ConstrID}. This would lead to additional restrictions on the form of the energy, entropy and their fluxes, but yield no residual entropy inequality. In particular, the entropy inequality on solutions \eqref{eq:gen:EntrOnSol} would vanish identically: $Q(Z)=0$. This analysis is performed in the following Section \ref{sec:eg2:adiab}.

%
%\medskip

%
%--------------------------------------------------------------------------------------
\subsection{Constitutive dependencies in the globally adiabatic case}\label{sec:eg2:adiab}

We now use the solution set approach to find the constraints on constitutive functions \eqref{eq:simplefluid:constitF:our:algo}
that will guarantee that the simple heat-conducting compressible anisotropic fluid \eqref{eq:fluid:gen}, \eqref{eq:fluid:2d:heatcond} undergoes an adiabatic process, that is, the entropy production inequality on all solutions (formula \eqref{eq:gen:EntrOnSol}) will be identically zero (cf. \eqref{eq:entropy_ineq:Liufluid:glob_adiab}). This is achieved by extending the set of constraint equations  \eqref{eq:OurAlgo:RunEx:ConstrID} with the left-hand side of the residual entropy inequality \eqref{eq:Liu:RunEx:ResIneq}. In symbolic computations, the requirement that this set of expressions must be zero is implemented, and the cases are split, with a command similar to \eqref{eq:eg2:casespl:plot}:
\beq\label{eq:eg2:casespl:plot:adiab}
\barr
\verb|Split_WRT_entropy_quantities:=DEtools[rifsimp](|\\
\verb|    [ Solution_Set_Constraints[], Residual_Entropy_Ineq,|\\
\verb|      Entropy_Ineqality_denom<>0,|\\
\verb|      diff(S(R, W),W)<>0, diff(E(R, W),W)<>0 ],|\\
\verb|    [ T11(R, W, Wx, Wy), T12(R, W, Wx, Wy), T22(R, W, Wx, Wy),|\\
\verb|      Q1(R, W, Wx, Wy), Q2(R, W, Wx, Wy),|\\
\verb|      Phi1(R, W, Wx, Wy), Phi2(R, W, Wx, Wy) ],|\\
\verb|    mindim=1, casesplit);|\\[2ex]
\verb|DEtools[caseplot](Split_WRT_entropy_quantities, pivots);|
\earr
\eeq
The last command is used to plot the case tree. As usual, the classifying functions are the constitutive functions not included in the list of arguments of \verb|rifsimp|; we again classify with respect to the forms of the entropy and energy density, $\eta(\rho, \vartheta)$ and $\epsilon(\rho, \vartheta)$.

\begin{figure}[!h]
  \begin{center}
  \includegraphics[width=0.4\textwidth]{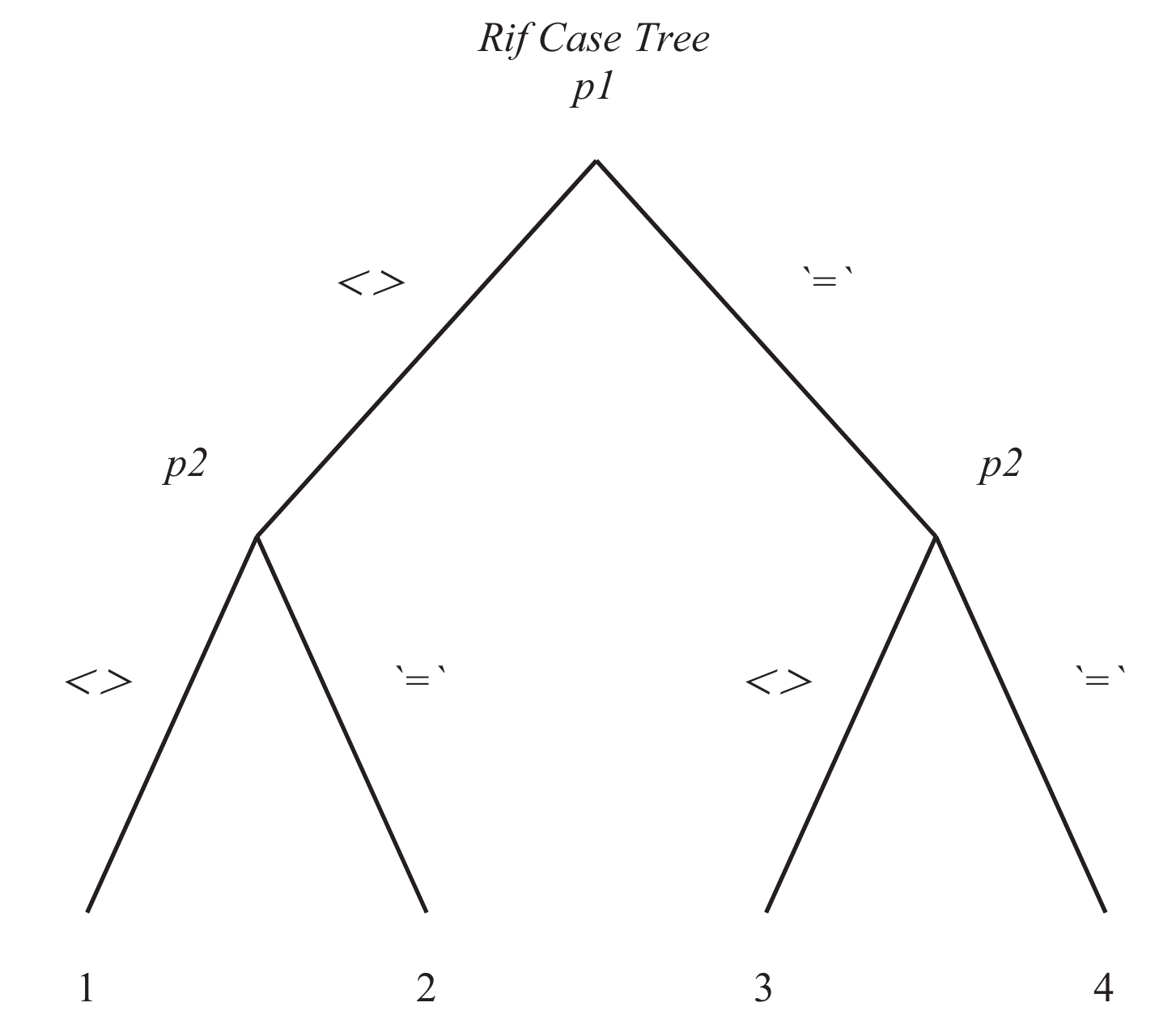}\\
  \end{center}
  \caption{Case tree for the classification of restrictions on the constitutive functions of the fully adiabatic two-dimensional heat-conducting fluid, Section \ref{sec:eg2:adiab}. The pivots are given by \eqref{eq:eg2:pivots}.}\label{fig:eg2:adiab}
\end{figure}

The fully adiabatic requirement, where \verb|Residual_Entropy_Ineq| must vanish in addition to \verb|Solution_Set_Constraints[]|,  naturally leads to more strict constraint equations than the general non-adiabatic case of Section \ref{sec:eg2:comp}. The new case tree is shown in Figure \ref{fig:eg2:adiab}. It involves four cases, determined by two pivots
\beq\label{eq:eg2:pivots}
\verb|p1|=\dfrac{\partial \epsilon}{\partial \vartheta} \dfrac{\partial^2 \eta}{\partial \rho\partial \vartheta} - \dfrac{\partial \eta}{\partial \vartheta} \dfrac{\partial^2 \epsilon}{\partial \rho\partial \vartheta}\,
,\qquad \verb|p2|=\dfrac{\partial \epsilon}{\partial \vartheta} \dfrac{\partial^2 \eta}{\partial \vartheta^2} - \dfrac{\partial \eta}{\partial \vartheta} \dfrac{\partial^2 \epsilon}{\partial \vartheta^2}.
\eeq

For the four cases that arise, one still obtains the same isotropy condition and the hydrostatic pressure as per \eqref{eq:eg2:isotr}.

\medskip\noindent\textbf{Case 1.}  In the general Case 1, with none of the pivots \eqref{eq:eg2:pivots} being equal to zero, the forms of the entropy and energy density $\epsilon(\rho,\vartheta)$, $\eta(\rho,\vartheta)$ are arbitrary, whereas the entropy and energy flux components are restricted by the conditions \eqref{eq:eg2:isotr}--\eqref{eq:eg2:gen:Qs} and  additional conditions
%\beq\label{eq:eg2:case1:fluxes} --- These are SAME as for non-adiab gen isotr fluid case
%\barr
%\dfrac{\partial q_i}{\partial (\partial_i \vartheta)}=\dfrac{\partial \Phi_i}{\partial (\partial_i \vartheta)}=0, \qquad \dfrac{\partial q_i}{\partial \rho} = \dfrac{{\partial \epsilon}/{\partial \vartheta}}{{\partial \eta}/{\partial \vartheta}}\,\dfrac{\partial \Phi_i}{\partial \rho}, \qquad i=1,2,\\[3ex]
%\dfrac{\partial q_2}{\partial (\partial_1 \vartheta)}+\dfrac{\partial q_1}{\partial (\partial_2 \vartheta)}=0,\qquad \dfrac{\partial \Phi_2}{\partial (\partial_1 \vartheta)}+\dfrac{\partial \Phi_1}{\partial (\partial_2 \vartheta)}=0,\\[3ex]
%\earr
%\eeq
\beq\label{eq:eg2:case1:fluxes:A}
\barr\partial_i \vartheta\,\dfrac{\partial q_i}{\partial \vartheta} \dfrac{\partial \eta}{\partial \vartheta} =
\partial_j \vartheta\,\dfrac{\partial \Phi_j}{\partial \vartheta} \dfrac{\partial \epsilon}{\partial \vartheta},
\earr
\eeq
\beq\label{eq:eg2:case1:fluxes:B}
\barr
\partial_i \vartheta\,\dfrac{\partial \Phi_i}{\partial \rho}\verb|p2| = \partial_j \vartheta\,\dfrac{\partial \Phi_j}{\partial \vartheta}\verb|p1|,
\earr
\eeq
where summations in repeated indices are assumed.

\medskip\noindent\textbf{Case 2.}  Here, the entropy and energy densities are related by $\verb|p2|=0$, which yields
\beq\label{eq:eg2:case2:es}
      \epsilon(\rho,\vartheta) = A(\rho) + B(\rho)\, \eta(\rho,\vartheta),
\eeq
where $A(\rho)$, $B(\rho)$ and the entropy density $\eta(\rho,\vartheta)$ are arbitrary functions. With this, one has simpler and less restrictive  conditions on the entropy fluxes compared to those in Case 1; namely, instead of \eqref{eq:eg2:case1:fluxes:B}, they must satisfy
\[
\partial_j \vartheta\,\dfrac{\partial \Phi_j}{\partial \vartheta}=0,
\]
where summation is assumed.

\medskip\noindent\textbf{Case 3.} For this case, the relation between the entropy and energy densities is given by $\verb|p1|=0$, which implies
\beq\label{eq:eg2:case3:es}
      \epsilon(\rho,\vartheta) = A(\rho) + \int B(\vartheta) \dfrac{\partial \eta}{\partial \vartheta} \, d\vartheta,
\eeq
where $A(\rho)$, $B(\vartheta)$ and the entropy density $\eta(\rho,\vartheta)$  are arbitrary functions. In this case, the conditions on the fluxes are significantly less restrictive than the general adiabatic conditions \eqref{eq:eg2:Phis}, \eqref{eq:eg2:gen:Qs}, \eqref{eq:eg2:case1:fluxes:A}, \eqref{eq:eg2:case1:fluxes:B} of Case 1. In particular, second derivatives of the flux components  of opposite gradient components do not need to vanish, but should satisfy
\beq\label{eq:eg2:case3:fluxes}
\barr
\dfrac{\partial \eta}{\partial \vartheta}\dfrac{\partial^2 q_1}{\partial (\partial_2 \vartheta)^2}=\dfrac{\partial \epsilon}{\partial \vartheta}\dfrac{\partial^2 \Phi_1}{\partial (\partial_2 \vartheta)^2},\qquad \dfrac{\partial \eta}{\partial \vartheta}\dfrac{\partial^2 q_2}{\partial (\partial_1 \vartheta)^2}=\dfrac{\partial \epsilon}{\partial \vartheta}\dfrac{\partial^2 \Phi_2}{\partial (\partial_1 \vartheta)^2}.
\earr
\eeq
The dependence of the energy and entropy fluxes on $\rho$ is provided by he equations
\beq\label{eq:eg2:case3:fluxes:A}
\barr
\dfrac{\partial \eta}{\partial \vartheta}\,\dfrac{\partial q_i}{\partial \rho} = \dfrac{\partial \epsilon}{\partial \vartheta}\, \dfrac{\partial \Phi_i}{\partial \rho},\qquad i=1,2,\\[3ex]
\earr
\eeq
and the vanishing sum condition
\[
\partial_i \vartheta\,\dfrac{\partial \Phi_i}{\partial \rho}=0.
\]
Moreover, \eqref{eq:eg2:case1:fluxes:A} holds, and fluxes also satisfy the double summation relation
\beq\label{eq:eg2:case3:fluxes:B}
\barr
\partial_i \vartheta\, \partial_j \vartheta \, \dfrac{\partial \Phi_i}{\partial (\partial_j \vartheta)} =0.
\earr
\eeq

\medskip\noindent\textbf{Case 4.}  Finally, in the fourth case, the relation between the entropy and energy densities is given by $\verb|p1|=\verb|p2|=0$, which yields
\beq\label{eq:eg2:case4:es}
      \epsilon(\rho,\vartheta) = A(\rho) + C_1 \eta(\rho,\vartheta),
\eeq
where $A(\rho)$ and the entropy density $\eta(\rho,\vartheta)$ are arbitrary, and $C_1=\const$. In this case, the constraint equations can be completely solved; they yield
\beq\label{eq:eg2:case4:fluxes}
q_1(\rho, \vartheta, \partial_i \vartheta) = C_2 \Phi_1- G(\vartheta)\partial_2\vartheta,\qquad q_2(\rho, \vartheta, \partial_i \vartheta) = C_2 \Phi_2+ G(\vartheta)\partial_1\vartheta.
\eeq
Furthermore, $G(\vartheta)$ and $\Phi_{i}=\Phi_{i}(\rho, \vartheta, \partial_i \vartheta)$ are arbitrary functions, $C_2$ is an arbitrary constant. The stress tensor components are still given by \eqref{eq:eg2:isotr}; in particular, this case is \emph{barotropic}, i.e. the pressure is a function of density only:
\beq\label{eq:eg2:case4:P}
P=P(\rho)=\rho^2\,\dfrac{d A(\rho)}{d\rho}.
\eeq

\begin{remark}
The four cases considered above are complementary -- any one of the four cases
yields a different set of restrictions on the constitutive functions. As usual, fewer restrictions on the classifying functions $\epsilon(\rho,\vartheta)$, $\eta(\rho,\vartheta)$ correspond to more significant restrictions in the remaining constitutive functions (flux components), and vice versa. Families of constitutive functions satisfying the constraint equations of any of the cases are consistent with the zero entropy production requirement $Q(Z)= 0$ on solutions, see \eqref{eq:gen:EntrOnSol}.

\end{remark}

%-------------------------------------------------------------------------------------------------------------------------------------------------
%-------------------------------------------------------------------------------------------------------------------------------------------------
\section{Computational example 3: Two-dimensional non-simple fluid}\label{sec:Maple:EG3}

We now present an example in which differential consequences become important: a \emph{non-simple two-dimensional heat-conducting fluid}, given by the PDE system
\begin{subequations}\label{eq:fluid:2d:heatcond:nonsimp}
%\beq\label{eq:fluid:2d:heatcond:nonsimp}
%\barr
\beq\label{eq:fluid:2d:heatcond:nonsimp:rho}
\Pi^{\rho}:~\partial_t \rho  +\partial_1 (\rho u)+ \partial_2(\rho v)=0,
\eeq
\beq\label{eq:fluid:2d:heatcond:nonsimp:u}
\Pi^{u}:~\rho (\partial_t u +u \partial_1 u + v \partial_2 u) - \partial_1 T_{11} - \partial_2 T_{12} =0,
\eeq
\beq\label{eq:fluid:2d:heatcond:nonsimp:v}
\Pi^{v}:~\rho (\partial_t v +u \partial_1 v + v \partial_2 v) - \partial_1 T_{12} - \partial_2 T_{22} =0,
\eeq
\beq\label{eq:fluid:2d:heatcond:nonsimp:e}
\Pi^{\epsilon}:~\rho (\partial_t \epsilon + u \partial_1 \epsilon+ v \partial_2 \epsilon ) +  \partial_1 q_1 +  \partial_2 q_2 - T_{11} \partial_1 u - T_{12} \partial_2 u- T_{12} \partial_1 v- T_{22} \partial_2 v=0,
\eeq
\end{subequations}
where $\rho(t,x,y)$ is the fluid density, and $u(t,x,y)$, $v(t,x,y)$ are the velocity components, $T_{ij}$ are the components of the stress tensor, $q_i$ are the heat flux components, $\epsilon(t,x,y)$ is the internal energy density, and $\vartheta(t,x,y)$ is the fluid temperature. The entropy inequality is given by
\begin{equation}\label{eq:entropy_ineq:eg3:nonsimp}
\Pi^{\eta}:~\rho (\partial_t \eta + u \partial_1 \eta +v \partial_2 \eta) +  \partial_1 \Phi_1 +  \partial_2 \Phi_2 \geq 0,
\end{equation}
The above model is similar to that of Section \ref{sec:eg2}, the difference being the ``non-simple''  constitutive dependencies that involve time derivatives:
\beq\label{eq:eg3a:constit:dep}
\barr
\epsilon=\epsilon(\rho,\partial_t \rho, \vartheta),\quad \eta=\eta(\rho,\partial_t \rho, \vartheta),\\[2ex]
T_{ij}=T_{ij}(\rho, \vartheta), \quad q_i=q_i(\rho, \vartheta), \quad \Phi_{i}=\Phi_{i}(\rho, \vartheta).
\earr
\eeq
These constitutive dependencies will require the use of differential consequences to set up and analyze the entropy inequality. In the current example, the solution set algorithm behaves consistently, yielding a set of constraint equations that guarantee the nonnegative entropy production, while the M\"{u}ller-Liu procedure yields an overdetermined, overly restrictive set of Liu identities. This is the case because the M\"{u}ller-Liu procedure is based on linear algebra techniques, and does not take into account nonlinear terms and relationships (through differential consequences) between the higher derivatives on the solution set of the PDEs.

%
%In the second sub-example, we consider more complex dependencies of constitutive functions, concentrating on the cases when the fluid can be anisotropic: $T_{12}\ne 0$. Here we show that the M\"{u}ller-Liu procedure yields an underdetermined set of Liu identities. As a result, the constraint equations that follow do not guarantee the nonnegative entropy production -- it may not hold on \emph{all} solutions of the given model. Conversely, the solution set algorithm yields a set of \emph{sufficient conditions} for the nonnegative entropy production postulate to hold for \emph{all} solutions of all admissible models satisfying those conditions.

In this section, we do not list all the commands for the symbolic computations; they are largely parallel to those in the previous examples. For simplicity of presentation, where appropriate, we freely use the coordinate notation $x= x_1$, $y= x_2$, partial derivatives $\partial_x= \partial_1$, $\partial_y= \partial_2$, the velocity components $u= v_1$, $v= v_2$, and summation in repeated indices.

%\subsection{Example 3a: the simplest dependence on temperature gradient.}

\subsection{The solution set method}
%\medskip\noindent \textbf{(i) The solution set method.}

In \verb|Maple|/\verb|GeM| symbolic software, the equations \eqref{eq:fluid:2d:heatcond:nonsimp} and the variables they involve are declared in the symbolic software, using commands similar to those in Steps A, B, C of Section \eqref{sec:eg2}.

First, the model PDEs \eqref{eq:fluid:2d:heatcond:nonsimp} with the constitutive dependencies \eqref{eq:eg3a:constit:dep} are converted to the symbolic form using commands similar to \eqref{eg2:symb:setup}--\eqref{eq:get:Pi:1234}, and are stored in the variables
\[
\verb|Pi_R_Symb|,~~\verb|Pi_U_Symb|,~~\verb|Pi_V_Symb|,~~\verb|Pi_E_Symb|.
\]

Second, we need to determine an appropriate set of leading derivatives $\phi_L$ for this model, such that the model equations (and their differential consequences) can be solved for elements of $\phi_L$ (and their differential consequences). It can be shown that the leading derivatives can be chosen to be the time derivatives of the four scalar fields:
\beq \label{eq:eg3a:LD}
   \phi_L = \{\partial_t \rho, \partial_t u,  \partial_t v,  \partial_t \theta\}.
\eeq
Indeed, $\partial_t \rho$ can be found from $\Pi^{\rho}$ \eqref{eq:fluid:2d:heatcond:nonsimp:rho}
\beq \label{eq:eg3a:LD:Rt}
\partial_t \rho =-\partial_1 (\rho u) - \partial_2(\rho v),
\eeq
using the \verb|Maple| command
\[ % used: file Eg3_NonSimpleFluid_SolSetAppr_eg3a.mw !!!
\barr
\verb|Leading_ders_1:={Rt};|\\
\verb|solved_LD_1:=solve([Pi_R_Symb], Leading_ders_1);|\\
\earr
\]
and the solved form of the equation is thus stored in the variable \verb|solved_LD_1|. The velocity time derivatives $\partial_t u$, and $\partial_t v$ can be computed from  the PDEs $\Pi^{u}$ and $\Pi^{v}$ \eqref{eq:fluid:2d:heatcond:nonsimp:u}, \eqref{eq:fluid:2d:heatcond:nonsimp:v}
\beq \label{eq:eg3a:LD:UVtt}
\barr
\partial_t u  = -(u \partial_1 u + v \partial_2 u) - \rho^{-1}(\partial_1 T_{11} - \partial_2 T_{12}), \\[2ex]
\partial_t v  = -( u \partial_1 v + v \partial_2 v) - \rho^{-1}(\partial_1 T_{12} - \partial_2 T_{22}), \\[2ex]
\earr
\eeq
using the commands
\[ % used: file Eg3_NonSimpleFluid_SolSetAppr_eg3a.mw !!!
\barr
\verb|Leading_ders_2:={Ut,Vt};|\\
\verb|solved_LD_1:=solve([Pi_U_Symb, Pi_V_Symb], Leading_ders_2);|\\
\earr
\]

We observe that, due to the dependence of $\epsilon$ and $\eta$ on $\partial_t \rho$ in \eqref{eq:eg3a:constit:dep}, the energy equation $\Pi^{\epsilon}$ \eqref{eq:fluid:2d:heatcond:nonsimp:e} and the entropy condition $\Pi^{\eta}$ \eqref{eq:entropy_ineq:eg3:nonsimp} involve the second-order derivatives
\beq\label{eq:eg3a:lead:r2}
\partial^2_t \rho,\quad \partial_t\partial_1 \rho, \quad\partial_t\partial_2 \rho.
\eeq
The derivatives \eqref{eq:eg3a:lead:r2} are related with other higher-order derivatives of the physical fields; these relations are given by differential consequences of the field equations. In particular, it is necessary to use the corresponding differential consequences of the continuity equation $\Pi^{\rho}$ \eqref{eq:fluid:2d:dens}, namely, its derivatives by $t$, $x$, and $y$:
\begin{equation}\label{eq:eg3a:difcons:rho}
\barr
\partial_t\Pi^{\rho}:~~\partial^2_t \rho  +\partial_t \partial_i (\rho v_i)=0,\\[2ex]
\partial_1\Pi^{\rho}:~~\partial_t \partial_1 \rho  +\partial_1 \partial_i (\rho v_i)=0,\\[2ex]
\partial_2\Pi^{\rho}:~~\partial_t \partial_2 \rho  +\partial_2 \partial_i (\rho v_i)=0;
\earr
\end{equation}
these DEs are obtained in symbolic form using commands
\[
\barr
\verb|Pi_Rt_Symb:=gem_total_derivative(Pi_R_Symb,t);|\\
\verb|Pi_Rx_Symb:=gem_total_derivative(Pi_R_Symb,x);|\\
\verb|Pi_Ry_Symb:=gem_total_derivative(Pi_R_Symb,y);|
\earr
\]
The last two PDEs in \eqref{eq:eg3a:difcons:rho} can be solved for the differential consequences $\partial_t\partial_1 \rho$, $\partial_t\partial_2 \rho$ of the leading derivative $\partial_t\rho$:
\[
\barr
\verb|Leading_ders_3:={Rtx,Rty};|\\
\verb|solved_LD_3:=solve([Pi_Rx_Symb,Pi_Ry_Symb], Leading_ders_3);|
\earr
\]
Now, the first PDE in \eqref{eq:eg3a:difcons:rho} also involves the second derivatives
\beq\label{eq:eg3a:lead:uv2}
\partial_t \partial_1 v_1, ~~\partial_t \partial_2 v_2;
\eeq
these are related with other derivatives through the differential consequences of the PDEs $\Pi^{u}$ \eqref{eq:fluid:2d:heatcond:nonsimp:u} and  $\Pi^{v}$ \eqref{eq:fluid:2d:heatcond:nonsimp:v} with respect to  $x$ and $y$, respectively:
\begin{equation}\label{eq:eg3a:difcons:UV}
\barr
\partial_1\Pi^{u}:~~\partial_1 (\rho \sg{D}_t u - \partial_j  T_{1j})=0,\\[2ex]
\partial_2\Pi^{v}:~~\partial_2 (\rho \sg{D}_t v - \partial_j  T_{2j})=0;
\earr
\end{equation}
their symbolic form is computed using
\[
\barr
\verb|Pi_Ux_Symb:=gem_total_derivative(Pi_U_Symb,x);|\\
\verb|Pi_Vy_Symb:=gem_total_derivative(Pi_V_Symb,y);|
\earr
\]
The differential consequences \eqref{eq:eg3a:difcons:UV} of the momentum equations are now solved for the differential consequences $\partial_t\partial_1 u$, $\partial_t\partial_2 v$ of the leading derivatives $\partial_t u$,  $\partial_t v$. This is done after the substitution of previously solved for expressions $\partial_t\rho$, $\partial_t u$, and $\partial_t v$, $\partial_t\partial_1 \rho$, $\partial_t\partial_2 \rho$:
\[
\barr
\verb|Leading_ders_4:={Utx,Vty};|\\
\verb|solved_LD_4:=solve(subs(solved_LD_1 union solved_LD_2 union solved_LD_3, |\\
\verb|                        [Pi_Ux_Symb,Pi_Vy_Symb]) ,Leading_ders_4);|
\earr
\]
The resulting expressions for these are omitted due to their length. So far, one has expressions for $\partial_t\rho$, $\partial_t u$, and $\partial_t v$, $\partial_t\partial_1 \rho$, $\partial_t\partial_2 \rho$, $\partial_t\partial_1 u$, and $\partial_t\partial_2 v$ in terms of other quantities. Substituting these in the differential consequence $\partial_t\Pi^{\rho}$ of \eqref{eq:eg3a:difcons:rho}, one obtains an expression for $\partial^2_t \rho$:
\[
\barr
\verb|Leading_ders_5:={Rtt};|\\
\verb|solved_LD_5:=solve(|\\
\verb|    subs(solved_LD_1 union solved_LD_2 union solved_LD_3 union solved_LD_4 ,Pi_Rt_Symb), |\\
\verb|    Leading_ders_5);|
\earr
\]
%(We again omit the resulting expression because of its length and no immediate importance for the presentation.)

Finally, we have solved for the leading derivatives
\beq\label{eq:eg3a:LD15solved}
\partial_t \rho, \quad \partial_t u,  \quad \partial_t v
\eeq
of the set \eqref{eq:eg3a:LD}, and their differential consequences
\beq\label{eq:eg3a:LD15solved2}
\partial_t \partial_1 v_1, \quad \partial_t \partial_2 v_2, \quad \partial^2_t \rho,\quad \partial_t\partial_1 \rho, \quad\partial_t\partial_2 \rho,
\eeq
in terms of other variables and derivatives. The right-hand sides of the corresponding expressions do not contain these leading derivatives and differential consequences. The remaining leading derivative in \eqref{eq:eg3a:LD}, namely, $\partial_t \theta$, is found in terms of non-leading derivatives from the energy balance equation $\Pi^{\epsilon}$  \eqref{eq:fluid:2d:heatcond:nonsimp:e}, after the substitution of the expressions \verb|solved_LD_2|, \ldots, \verb|solved_LD_5| for the quantities \eqref{eq:eg3a:LD15solved}, \eqref{eq:eg3a:LD15solved2}:
\[
\barr
\verb| Leading_ders_6:={Wt};|\\
\verb| solved_LD_6:=solve(|\\
\verb|        subs(solved_LD_2 union solved_LD_3 union solved_LD_4 union solved_LD_5, Pi_E_Symb),|\\
\verb|        Leading_ders_6);|
\earr
\]
We note that the expression for $\partial_t \rho$ ,\verb|solved_LD_1| is not used, in order to not spoil the dependency of $\epsilon=\epsilon(\rho,\partial_t \rho, \vartheta)$ in \eqref{eq:eg3a:constit:dep}, but it must be used later. As a result, the solution set of the PDE model of a non-simple fluid \eqref{eq:fluid:2d:heatcond:nonsimp} is fully described, up to the required derivative order, by the \emph{solution set conditions}: relationships \verb|solved_LD_1|, \ldots, \verb|solved_LD_6| of the  leading derivatives \eqref{eq:eg3a:LD} and their differential consequences \eqref{eq:eg3a:LD15solved2}, together with other jet space variables of the given model.

\medskip

%with the right-hand sides not involving these derivatives. These expressions define the solution set of the current physical model of the two-dimensional non-simple fluid in the jet space.

Having established the relations defining the solution set as a manifold in the jet space,  we now derive the constraints on the constitutive functions that follow from the requirement of nonnegative entropy production. The symbolic version of the entropy inequality \eqref{eq:entropy_ineq:eg3:nonsimp} is computed using commands similar to \eqref{eq:eg2:PiS}. In order to restrict the entropy inequality to the solution set of the given model, one substitutes into it all the solution set conditions, except for the expression for $\partial_t \rho$, which will be done later, and is not done at this stage to prevent the spoiling of the constitutive dependency of the entropy density $\eta=\eta(\rho,\partial_t \rho, \vartheta)$. The following commands are used:
\[
\barr
\verb|Sol_Set_Conditons:= solved_LD_2 union solved_LD_3|\\
\verb|          union solved_LD_4 union solved_LD_5 union solved_LD_6;|
\earr
\]
\[
\verb|Entropy_Ineqality_Sol_Set:=simplify(subs(Sol_Set_Conditons, Pi_S_Symb));|
\]
The numerator and denominator of the entropy inequality are obtained like before by commands \eqref{eg1:step:E:3}. %on the solution set of the model
In particular, the denominator \verb|Entropy_Ineqality_denom| is given by $\partial \epsilon/\partial \vartheta$, and the numerator will be further analyzed.

The remaining essential step to ensure that the entropy inequality is analyzed on the solution set of the model is to use the continuity equation $\Pi^{\rho}$ \eqref{eq:fluid:2d:dens}, by substituting into entropy inequality numerator \verb|Entropy_Inequality_numer| the right-hand side of the relation \eqref{eq:eg3a:LD:Rt} instead of the leading derivative $\partial_t \rho$. After that, the constitutive dependence of  $\eta=\eta(\rho,\partial_t \rho, \vartheta)$,  $\epsilon=\epsilon(\rho,\partial_t \rho, \vartheta)$ that will be spoiled by these substitutions can be restored. In terms of \verb|Maple|-based symbolical computations, This is done using the following forward and backward substitutions.
 %Eg3_NonSimpleFluid_SolSetAppr_eg3a_Final.mw %\beq\label{eg3:entrIneq:numer:UseRt}
\[
\barr
\verb|Entropy_Inequality_numer_b:=expand(subs( solved_LD_1,|\\
\verb|                       convert(Entropy_Inequality_numer,D)));|\\
\verb|Entropy_Inequality_numer_final:=convert(subs(|\\
\verb|                       rhs(solved_LD_1[1])=lhs(solved_LD_1[1]),|\\
\verb|                       Entropy_Inequality_numer_b), diff);|\\
\earr
\]
The variable \verb|Entropy_Inequality_numer_final| now contains the numerator of the entropy inequality computed on the solution set of the model.

The problem includes the \emph{free elements}
\beq\label{eg3a:arb:el}
\barr
t, x, y,~~u, v,~~ \partial_i \rho, \partial_i u, \partial_i v, \partial_i \theta,\\[1ex]
\partial_i\partial_j \rho, \partial_i\partial_j u, \partial_i\partial_j v, \partial_i\partial_j \theta, \\[1ex]
\partial_t^2 u, \partial_t^2 v,  \partial_t^2 \theta,~~\partial_t\partial_i \theta,~~ \partial_t\partial_y u, \partial_t\partial_x v.
\earr
\eeq
which are all variables and their derivatives with the exception of the leading derivatives, their differential consequences, and the constitutive dependencies -- the free elements of \eqref{eg3a:arb:el} are computed using the command \eqref{eg1:step:E:2}, being placed into the variable \verb|Arbitrary_Elements|. On the solution set of the model, the free elements are indeed arbitrary.

The constraint equations are now obtained by setting to zero the coefficients at all combinations of free elements \eqref{eg3a:arb:el} in the entropy inequality numerator \verb|Entropy_Inequality_numer_final|:
\[
\barr
\verb|Solution_Set_Constraints:=convert({coeffs_constraints} minus {Residual_Entropy_Ineq},list);|
\earr
\]
The resulting 26 constraint equations are redundant, as is commonly the case, and can be written as
\beq\label{eq:3a:ConstrID}
\barr
%T_{12}\dfrac{\partial \eta}{\partial\vartheta} =0;\\[3ex] %1 in Eg3_NonSimpleFluid_SolSetAppr_eg3a_Final3.mw

\dfrac{\partial \epsilon}{\partial\vartheta} \dfrac{\partial \Phi_i}{\partial\rho} -\dfrac{\partial \eta}{\partial\vartheta} \dfrac{\partial q_i}{\partial\rho}=0,\quad \dfrac{\partial \epsilon}{\partial\vartheta} \dfrac{\partial \Phi_i}{\partial\vartheta} -\dfrac{\partial \eta}{\partial\vartheta} \dfrac{\partial q_i}{\partial\vartheta}=0, \quad i=1,2;\\[3ex] %2,3,4,5

\dfrac{\partial \epsilon}{\partial\vartheta} \dfrac{\partial \eta}{\partial(\partial_t\rho)} -\dfrac{\partial \eta}{\partial\vartheta} \dfrac{\partial \epsilon}{\partial(\partial_t\rho)}=0;\\[3ex] %6-9, 12-26

%T_{11}\dfrac{\partial \eta}{\partial\vartheta} - \rho^2\left(\dfrac{\partial \epsilon}{\partial\vartheta} \dfrac{\partial \eta}{\partial\rho} -\dfrac{\partial \eta}{\partial\vartheta} \dfrac{\partial \epsilon}{\partial\rho}\right)=0;\quad T_{22}\dfrac{\partial \eta}{\partial\vartheta} - \rho^2\left(\dfrac{\partial \epsilon}{\partial\vartheta} \dfrac{\partial \eta}{\partial\rho} -\dfrac{\partial \eta}{\partial\vartheta} \dfrac{\partial \epsilon}{\partial\rho}\right)=0;\\[3ex] %10, 11

\rho^2\left(\dfrac{\partial\epsilon}{\partial\rho}\dfrac{\partial \eta}{\partial\vartheta} - \dfrac{\partial \epsilon}{\partial\vartheta} \dfrac{\partial \eta}{\partial\rho} \right)  \delta_{ij}+\dfrac{\partial \eta}{\partial\vartheta} T_{ij}=0,\quad i,j=1,2;\\[3ex] %1, 10, 11

\partial  \epsilon/ \partial \vartheta\ne 0,
\earr
\eeq
where the last condition implies that the entropy inequality denominator \verb|Entropy_Ineqality_denom| does not vanish.

The residual entropy inequality corresponds to the ``free coefficient'' -- the part of the entropy inequality independent of the free elements \eqref{eg3a:arb:el}. It is computed by setting the free elements to zero, with commands
\[
\barr
\verb|Residual_Entropy_Ineq:=expand(simplify(eval(subs(|\\
\verb|            {map(x->x=0,Arbitrary_Elements)[]}, |\\
\verb|            Entropy_Inequality_numer_final))));|
\earr
\]
and is identically equal to zero.

The results obtained for the non-simple fluid model  \eqref{eq:fluid:2d:heatcond:nonsimp} with the constitutive dependencies \eqref{eq:eg3a:constit:dep} can be summarized as follows:
\begin{enumerate}

  \item The requirement that the entropy inequality is independent on all parametric (free) derivatives, for the current model and the posed constitutive dependencies, yields not only a sufficient conditions for the process to have nonnegative entropy production, but in fact an \emph{adiabatic process}. This is the case since the residual entropy inequality vanishes identically.

  \item Moreover, under the constitutive dependence assumptions \eqref{eq:eg3a:constit:dep}, the two-dimensional non-simple fluid must again be isotropic.

  \item The constraint equations \eqref{eq:3a:ConstrID} can be integrated, and lead to the admissible constitutive function forms
    \beq\label{eq:eg3:sol}
    \barr
        \epsilon(\rho,\partial_t \rho, \vartheta) = F(\rho)+C_0 \eta(\rho,\partial_t \rho, \vartheta);\\[2ex]
        q_i(\rho, \vartheta) = C_0\Phi_{i}(\rho, \vartheta) + C_i,\quad i=1,2;\\[2ex]
        T_{12}(\rho, \vartheta)=0,\\[2ex]
        T_{11}(\rho, \vartheta)=T_{22}(\rho, \vartheta)=-P(\rho),\\[2ex]
        P(\rho)=\rho^2 F'(\rho),
    \earr
    \eeq
    where $\eta(\rho,\partial_t \rho, \vartheta)$ and $\Phi_{i}(\rho, \vartheta)$ are arbitrary functions, and $C_0, C_i$ are arbitrary constants.

\end{enumerate}

It is evident that the results of the current section are different from those obtained for the simple heat-conducting compressible anisotropic fluid, see the running example of Section \ref{sec:EntropyPrinciple:ML}, due to different constitutive dependencies. We note that similarly to the simple fluid example, here again $T_{12}(\rho, \vartheta)=0$, hence the two-dimensional non-simple fluid must again isotropic.

\begin{remark} We note that case splitting/classification could also have been done for this example, but was omitted for brevity: only the most general case was considered and presented.
\end{remark}

%\medskip\noindent \textbf{(ii) Comparison with the classical M{\"u}ller-Liu procedure.}

%
%-----------------------------------------------------------------------------
\subsection{Comparison with the classical M{\"u}ller-Liu procedure}

If the classical M{\"u}ller-Liu procedure is applied to find the constraints on the  constitutive functions \eqref{eq:eg3a:constit:dep} for the non-simple two-dimensional heat-conducting fluid model, one chooses the Lagrange multiplier ansatz consistent with \eqref{eq:eg3a:constit:dep}:
\[
\Lambda^{\rho}=\Lambda^{\rho}(\rho,\partial_t \rho, \vartheta),\qquad\Lambda_i^{v}=\Lambda_i^{v}(\rho,\partial_t \rho, \vartheta), \qquad
\Lambda^{\epsilon}=\Lambda^{\epsilon}(\rho,\partial_t \rho, \vartheta).
\]
In the extended entropy inequality, the coefficients at the 35 ``arbitrary'' derivatives
\[
\partial_j \rho, ~~\partial_t\partial_j \rho,~~ \partial_j\partial_k \rho,~~
\partial_j v_i, ~~ \partial_t v_i,   ~~\partial_t\partial_k v_i, ~~\partial_j\partial_k v_i,~~
\partial_j \vartheta, ~~\partial_t\partial_j \vartheta,~~ \partial_j\partial_k \vartheta
\]
are set to zero. The resulting Liu identities do not take into account relationships between these derivatives on the solution set, since the differential consequences \eqref{eq:eg3a:difcons:rho}-\eqref{eq:eg3a:difcons:UV} are never used in the linear algebra-based Liu's lemma. As a result, the Liu identities are overdetermined (contain excessive restrictions), and at the same time may not provide a sufficient condition for nonnegative entropy production. In particular, they require constant fluxes
\beq\label{eq:eg3:Liu}
  q_i(\rho, \vartheta)=\const, ~~\Phi_{i}(\rho, \vartheta) =\const,
\eeq
which are unnecessary restrictions on the energy and entropy fluxes in comparison to the second equation of \eqref{eq:eg3:sol}.

%\newpage
%
%-------------------------------------------------------------------------------------------------------------------------------------------------
%-------------------------------------------------------------------------------------------------------------------------------------------------

%-------------------------------------------------------------------------------------------------------------------------------------------------
\section{Discussion and conclusions}\label{sec:discussion}

%***CLASSIF NOTE: fewer restrictions on the classifying functions correspond to more significant restrictions in the remaining constitutive functions (flux components), and vice versa. Seen in all classifications, for example, see computational example 2 (sec.4), specifically, its adiabatic classification. ***
%***NOTE: in last e.g., $T12<>0$: truly anisotropic

In the current work, a solution set-based entropy principle for constitutive modeling and its symbolic implementation was presented. It is applicable to any continuum model involving an entropy production-type inequality that has to hold for every solution of the model. The solution set approach generalizes the popular M{\"u}ller-Liu procedure to ensure the correct treatment of relationships of the higher-order derivatives of physical fields to the solution manifold of the model PDEs.

The solution set method is based on the solution of the dynamic field equations, as algebraic equations, for a set of chosen leading derivatives, in such a way that the right-hand sides of such solved expressions do not involve the leading derivatives or their differential consequences. If required, the field equations are solved for differential consequences of the leading derivatives, and appropriate substitutions are made. Similarly to what M{\"u}ller proposed in his later work (e.g., Ref.~\cite{muller1970new}), these expressions of leading derivatives (and possibly their differential consequences) in terms of all other quantities  are substituted into the entropy inequality, and coefficients at remaining free derivatives (the actual free elements) are required to vanish. So, similarly to the famous M{\"u}ller-Liu procedure, the solution set approach yields a set of constraint equations on the constitutive functions, and a residual entropy inequality. At the same time no artificial ``Lagrange multipliers'' are involved, and most importantly, the resulting constraint equations and the residual entropy inequality provide a set of \emph{necessary and sufficient conditions} for the nonnegative entropy production. We underline that this is generally not the case for the Liu identities following from the M{\"u}ller-Liu approach, see Table \ref{tab:principles}.

The solution set algorithm implementation in \verb|Maple| symbolic software, using the \verb|GeM| package, was presented and supported by several computational examples. In Section \ref{sec:Maple:EG1}, an elementary example based on one-dimensional gas dynamics is considered. We showed how the application of the solution set method leads to a classification problem that yields several possible cases of constitutive functions. To the best of our knowledge, such classifications have not appeared in constitutive modeling literature so far.

In Sections \ref{sec:eg2} and \ref{sec:Maple:EG3}, models of two-dimensional simple and non-simple heat-conducting fluids were considered. For the non-simple fluid, the constitutive dependencies involved a time derivative of the density, and thus, the entropy inequality required the use of a differential consequence in conjunction with the continuity equation. It was demonstrated that, while the solution set algorithm behaved consistently, returning a set of sufficient conditions for the nonnegative entropy production, the M\"{u}ller-Liu procedure led to an overly restrictive, mathematically inconsistent set of Liu identities.

%initially considered in Refs.~\cite{wang_comparison_1999,wang_shearing_1999} (see also
As a more computationally demanding example, a model of a two-dimensional motion of a granular solid was considered in Appendix \ref{sec:Maple:EG4}. In this example, we demonstrated a systematic computation of solution set conditions, and their use to obtain a set of constraint equations, including the corresponding symbolic computations. Due to the complex, extended constitutive dependencies, which included multiple fields, their derivatives, and symmetric stress tensor components, the resulting set of constraint equations does not yield any particularly simple relations between constitutive functions. Specifically, here it does not follow that the medium must be isotropic, as it was the case in Sections \ref{sec:eg2} and \ref{sec:Maple:EG3}; moreover, a simple relation \eqref{eq:phi:q:gran} between the flux components that was found in \cite{wang_comparison_1999} using the classical the  M{\"u}ller-Liu procedure was not found to be necessary for the nonnegative entropy production. Further analysis of the complex constraint equations obtained for this model is a direction of future work.

The important features of the new proposed solution set method for constitutive modeling and its symbolic implementation can be summarized as follows.
\begin{itemize}
  \item The constraint equations and the  residual entropy inequality following from the solution set method provide necessary and sufficient conditions for the nonnegative entropy production, unlike the Liu identities following from the M{\"u}ller-Liu approach.
  \item The solution set approach does not require Lagrange multipliers, and can, in principle, be applied to physical models and constitutive dependencies of arbitrary complexity.
  \item The symbolic implementation of the solution set method uses efficient mathematical techniques of differential simplification and elimination of constraint equations that are built into \verb|Maple| software, and, as illustrated in Sections  \ref{sec:Maple:EG1} and \ref{sec:eg2},  allows for automated classification and case splitting that lead to the determination of non-overlapping sets of admissible constitutive functions.
\end{itemize}

In conclusion, the proposed solution set method provides a mathematically consistent generalization for the well-known M{\"u}ller-Liu approach, is often computationally simpler, and therefore is a natural method of choice for constitutive modeling applications.

Possible future work directions include further computational examples of the solution set entropy principle applied to models of practical interest, both in dynamical setting and in thermodynamic equilibrium. Thermodynamic equilibrium is defined as a state or process in which both temperature and velocity are equally distributed and thus the entropy production is equal to zero (see, e.g., \cite[p. 198f]{hutter_continuum_2004}). With this, a set of additional constraints can be obtained, allowing for the deduction of the forms of so-called equilibrium parts of constitutive functions, that then have to be amended by often postulated dynamic parts (see, e.g., \cite{hess2017thermodynamically,schneider2009solid}). Moreover, since entropy principles often provide only basic restrictions on the constitutive functions, it is also of interest to use entropy principles in conjunction with other modeling approaches and to take into consideration the capabilities provided by computational tools. It is also of interest to study possible connections of the proposed methodology with the description of dynamics of internal variables in irreversible processes.

%---------------------------------------------------------------------------------------------------------------------------------------------------------------------------------------------
\subsubsection*{Acknowledgements}

The authors thank Professors Yongqi Wang and Martin Oberlack for ideas and discussions related to this project, and are also grateful for the financial support through the NSERC Discovery grant RGPIN-2014-05733, and the Deutsche Forschungsgemeinschaft (DFG) project WA 2610/3-1.

\bibliographystyle{ieeetr}
\footnotesize
\bibliography{Bib_RevisionF}
\normalsize

%-------------------------------------------------------------------------------------------------------------------

\begin{appendix}

\section{Solution set entropy principle for two-dimensional granular flows}\label{sec:Maple:EG4}

We now illustrate the application of the solution set approach to the analysis of the constitutive dependencies of a two-dimensional granular solid motion model, similar to the model proposed by Wang and Hutter \cite{wang_comparison_1999,wang_shearing_1999} (see also \cite{hess2017thermodynamically,schneider2009solid}). Compared to the usual fluid flow equations considered above, the granular flow model incorporates a supplementary balance equation for the solid volume fraction $\xi$, which also appears in the remaining governing equations. The physical fields of the model are given by
\[
\phi= \{\rho, v_i, \xi, \vartheta\},
\]
where $\rho$ is the true mass density, $v_1=u$ and $v_2=v$ are the velocity components, and $\vartheta$ is the temperature . The governing equations for mass, momentum, equilibrated forces and energy are given by
\begin{subequations}\label{eq:granular:eqs}
\begin{equation}\label{eq:granular:2d:dens}
\Pi^{\rho}:~ \partial_t ( \xi \rho) +  \partial_j (\xi \rho v_j)=0,
%\Pi^{v}: \rho D_t v - (\div T) =0,\qquad i=1\ldots, n,\\
\end{equation}
\begin{equation}\label{eq:granular:2d:mom}
%\Pi_i^{v}:~\xi  \rho  \partial_t  v_i +\xi  \rho v_j  \partial_j  v_i  - \partial_j  T_{ij}  =0,\quad i=1,2,
\Pi_i^{v}:~\xi  \rho  \, \sg{D}_t v_i  - \partial_j  T_{ij}  =0,\quad i=1,2,
\end{equation}
\begin{equation}\label{eq:granular:2d:nu}
\Pi^{\xi}:~\xi \rho\, \sg{D}_t^2 \xi  - \partial_j h_j -\xi \rho f=0,
\end{equation}
\begin{equation}\label{eq:fluid:2d:en}
\Pi^{\epsilon}:~\xi \rho \, \sg{D}_t\epsilon  - T_{jk} \partial_j v_k
- h_j  \partial_j ( \sg{D}_t\xi ) + \xi \rho f \,\sg{D}_t\xi + \partial_j q_j=0,
\end{equation}
\end{subequations}
where  $h_i$ are the equilibrated stress vector components, and $f$ is the equilibrated internal body force. The summations in repeated indices ($j$ and $k$, from 1 to 2) are assumed, as usual, in all PDEs where these indices occur, as well as in the material derivatives; we also use $x=x_1, y=x_2$, $\partial_1=\partial_x$, $\partial_2=\partial_y$. The second material derivative is given by $\sg{D}_t^2=(\partial_t + v_j  \partial_j)(\partial_t + v_k  \partial_k)$. The constitutive functions of the model \eqref{eq:granular:eqs} are given by
\beq\label{eq:gran:constitF}
\psi = \left\{ \epsilon, \eta, \Phi_i, T_{ij}, q_i, f, h_i\right\},\quad i=1,2
\eeq
where we note that in Ref.~\cite{wang_comparison_1999,wang_shearing_1999}, instead of the internal energy $\epsilon$, the free Helmholtz energy $\Psi_h=\epsilon -\vartheta \eta$
was taken to be a constitutive function.

The balance of equilibrated forces $\Pi^{\xi}$ was introduced by Goodman and Cowin \cite{goodman1972continuum} to describe the evolution of the volume fraction and, with this, the microstructure of the granular material in its development. Its introduction also yields additional contributions to the energy balance \eqref{eq:fluid:2d:en}.  The balance of angular momentum, as usual, yields the symmetry of the stress tensor, $T_{ij}=T_{j\,i}$. We do not append this requirement to equations, but rather treat the three independent components $T_{11}$, $T_{12}=T_{21}$, and $T_{22}$ accordingly.

Following Refs.~\cite{wang_comparison_1999,wang_shearing_1999}, let us assume a common constitutive dependency for all functions except the internal energy and entropy density. The latter two will have a reduced dependence, which does not include $\partial_t \xi$. With this assumption, the constitutive dependencies are given by
\beq\label{eq:granular2d:constit}
\barr
  f=   f (\xi, \partial_i \xi, \partial_t \xi, \rho, \vartheta, \partial_i \vartheta, D_{ij}),&\\[1ex]
   T_{ij}=   T_{ij}(\xi, \partial_i \xi, \partial_t \xi, \rho, \vartheta, \partial_i \vartheta, D_{ij} ), & q_i=   q_i(\xi, \partial_i \xi, \partial_t \xi, \rho, \vartheta, \partial_i \vartheta, D_{ij}),\\[1ex]
\Phi_{i}=  \Phi_{i}(\xi, \partial_i \xi, \partial_t \xi, \rho, \vartheta, \partial_i \vartheta, D_{ij}), & h_i=   h_i(\xi, \partial_i \xi, \partial_t \xi, \rho, \vartheta, \partial_i \vartheta, D_{ij}),\\[1ex]
\epsilon=  \epsilon(\xi, \partial_i \xi, \rho, \vartheta, \partial_i \vartheta, D_{ij}),& \eta=   \eta(\xi, \partial_i \xi, \rho, \vartheta, \partial_i \vartheta, D_{ij}),
\earr
\eeq
which includes the dependence on the components of the symmetric strain tensor
\beq\label{eq:gran:constitDep:Dij}
D_{ij}= \frac{1}{2} \left(\partial_j v_i + \partial_i v_j \right).
\eeq
%%%%%%%%%%%%%%%%%%

The entropy inequality for the current two-dimensional model is given by
\begin{equation}\label{eq:entropy_ineq:eg:Gran}
\Pi^{\eta}:~\Pi^{\eta}:~ \xi \rho \,\sg{D}_t \eta +  \partial_i \Phi_i \geq 0
\end{equation}
and is slightly different from that in previous examples due to the presence of the solid volume fraction.

%the same expression \eqref{eq:entropy_ineq:eg3:nonsimp} as before.

\subsection{The solution set entropy principle for granular flows}

%\smallskip\noindent\textbf{The solution set entropy principle.}

The above granular flow model is rather complex, in particular, due to extended constitutive dependencies \eqref{eq:granular2d:constit}. We now illustrate how one can obtain equations that fully describe the solution set of the problem including the necessary differential consequences and symmetry conditions, and then the steps required to further analyze the granular flow model in the context of the solution set entropy principle.

One may follow the same steps as in Section \ref{sec:Maple:EG1} and other examples; we present \verb|Maple| code where needed. The symbolic package initialization is the same as before, in \ref{eq:gem:initeq}. Next, the constitutive dependencies \eqref{eq:granular2d:constit} are declared, for example, by the following commands:
\[
\barr
\verb|ind:=t,x,y;   dep:=R(ind), U(ind), V(ind), W(ind), N(ind);|\\
\verb|Constit_Dependence:=N(ind),diff(N(ind),t),diff(N(ind),x),diff(N(ind),y),R(ind),|\\
\verb|                    W(ind),diff(W(ind),x),diff(W(ind),y),diff(U(ind),x),|\\
\verb|                    diff(U(ind),y),diff(V(ind),x),diff(V(ind),y);|\\
\verb|Constit_Dependence2:=N(ind),diff(N(ind),x),diff(N(ind),y),R(ind),|\\
\verb|                    W(ind),diff(W(ind),x),diff(W(ind),y),diff(U(ind),x),|\\
\verb|                    diff(U(ind),y),diff(V(ind),x),diff(V(ind),y);|\\
\verb|Constit_F:=T11(Constit_Dependence), T12(Constit_Dependence),|\\
\verb|           T22(Constit_Dependence), Q1(Constit_Dependence), |\\
\verb|           Q2(Constit_Dependence), Phi1(Constit_Dependence),|\\
\verb|           Phi2(Constit_Dependence), E(Constit_Dependence2),|\\
\verb|           H1(Constit_Dependence), H2(Constit_Dependence),|\\
\verb|           F(Constit_Dependence), S(Constit_Dependence2);|\\
\earr
\]
\[
\barr
\verb|gem_decl_vars(indeps=[ind], deps=[dep], freefunc=[Constit_F]);|
\earr
\]
where the variable \verb|Constit_Dependence2| corresponds to the reduced constitutive dependence of energy and entropy densities. The symbolic variables are denoted as usual, $\rho=\verb|R|$, $v_1=\verb|U|$, $v_1=\verb|U|$, $v_2=\verb|V|$, $\vartheta =\verb|W|$, and the solid volume fraction $\xi=\verb|N|$.

In symbolic computations, it is not convenient to consider functional dependencies on objects like the symmetrized derivatives \eqref{eq:gran:constitDep:Dij}. Instead all constitutive functions may, for now, be assumed to depend on \emph{all} velocity gradient components $\partial_j v_i$, and symmetrization conditions will be appended later.

The next step is to declare the model equations. Similarly to the previous examples, material derivative operators can be used to input the five scalar PDEs \eqref{eq:granular:eqs} in the symbolic form, naming them, for example,
\[
\verb|Pi_R_Symb|,~~\verb|Pi_U_Symb|,~~\verb|Pi_V_Symb|,~~\verb|Pi_N_Symb|,~~\verb|Pi_E_Symb|
\]
(see \eqref{eq:sec31:defeqs}, \eqref{eq:sec31:defeqs2} in Section \ref{sec:Maple:EG1:framew} and similar commands in Section \ref{sec:eg2:comp}).

On the next step, one computes the conditions that determine the solution set of the model. As the leading derivatives, it is natural to take the following (highest in time $t$ and in differential order) derivatives of all fields:
\beq \label{eq:eg4gran:LD}
   \phi_L = \{\partial_t \rho, \partial_t u,  \partial_t v,  \partial_t\partial_x \theta,  \partial_t^2 \xi\}.
\eeq
Then the model PDEs \eqref{eq:granular:eqs} can be solved for these quantities in a straightforward way,
\[
\barr
\verb|Leading_ders:={Rt,Ut,Vt,Ntt,Wtx};|\\
\verb|solved_DEs:=solve([Eq_R_Symb,Eq_U_Symb,Eq_V_Symb,Eq_N_Symb,Eq_E_Symb], Leading_ders);|
\earr
\]
The right-hand sides of the resulting expressions for $\partial_t \rho$, $\partial_t u$, $\partial_t v$, $\partial_t^2 \xi$ in \verb|solved_DEs| do not involve leading derivatives $\phi_L$ or their differential consequences. However, the expression for $\partial_t\partial_x \theta$ does involve the differential consequences $\partial_t\partial_x u$, $\partial_t \partial_x v$, $\partial_t\partial_y u$, $\partial_t \partial_y v$, due to the constitutive dependency of $\epsilon$. In a way parallel to what was done in Section \ref{sec:eg2}, one has to compute differential consequences of the momentum equations \eqref{eq:granular:2d:mom}:
\beq\label{eq:eg4gran:LD:difcons}
\partial_j\Pi_i^{v},\quad i,j=1,2,
\eeq
solve them for the differential consequences $\partial_t\partial_x u$, $\partial_t \partial_x v$, $\partial_t\partial_y u$, $\partial_t \partial_y v$ of leading derivatives \eqref{eq:eg4gran:LD}, and substitute this into the expression for $\partial_t\partial_x \theta$. The differential consequences \eqref{eq:eg4gran:LD:difcons} are computed as follows:
\[
\barr
\verb|Pi_Ux_Symb:=gem_total_derivative(Eq_U_An,x):|\\
\verb|Pi_Uy_Symb:=gem_total_derivative(Eq_U_An,y):|\\
\verb|Pi_Vx_Symb:=gem_total_derivative(Eq_V_An,x):|\\
\verb|Pi_Vy_Symb:=gem_total_derivative(Eq_V_An,y):|\\
\verb|difcons_subs:=simplify(subs(solved_DEs,{Pi_Ux_Symb,Pi_Uy_Symb,Pi_Vx_Symb,Pi_Vy_Symb})):|
%former DC_UV_2
\earr
\]
The variable \verb|difcons_subs| now contains the differential consequences \eqref{eq:eg4gran:LD:difcons} of the momentum equations, with leading derivatives \verb|Leading_ders| substituted. Now the differential consequences are solved,
\[
\verb|difcons_solved:=solve(difcons_subs, {Utx, Uty, Vtx, Vty}):|
\]
and substituted into the energy equation (the last one in \verb|solved_DEs|):
\[
\barr
\verb|energy_equation_solved:=solved_DEs[-1];|\\
\verb|energy_equation_solved_final:=simplify(eval(subs(solved_DEs_UVij,solved_DEs_Last))):|
\earr
\]
Finally, the solution set required for the entropy inequality \eqref{eq:entropy_ineq:eg:Gran} in conjunction with the constitutive dependencies \eqref{eq:granular2d:constit} is defined by the first four equations \verb|solved_DEs|, the energy equation solved for $\partial_t\partial_x \theta$, and differential consequences \eqref{eq:eg4gran:LD:difcons} of the momentum equations:
\beq\label{eq:eg4gran:all_sol_set_eqs_with_DC}
\barr
\verb|solution_set_equations_final:=|\\
\verb|   [solved_DEs[1..4][], energy_equation_solved_final, difcons_solved[] ]:|
\earr
\eeq

As a result, the symbolic vector variable \verb|solution_set_equations_final| \eqref{eq:eg4gran:all_sol_set_eqs_with_DC} contains the PDEs \eqref{eq:granular:eqs}, \eqref{eq:eg4gran:LD:difcons} solved for the prolonged set of leading derivatives
\beq \label{eq:eg4gran:LD:prol}
   \widetilde{\phi}_L = \{\partial_t \rho, \partial_t u,  \partial_t v,  \partial_t\partial_x \theta,  \partial_t^2 \xi, \partial_t\partial_x u, \partial_t \partial_x v, \partial_t\partial_y u, \partial_t \partial_y v\}
\eeq
%This is the file Eg4_Granular_SolSetAppr_07_appendix_Wtx.mw
in terms of all other quantities.

\begin{remark}
If for the granular flow model, other constitutive dependencies than \eqref{eq:granular2d:constit} are chosen, for example, if the internal energy and entropy density $\epsilon$, $\eta$ were independent of the temperature gradient $\partial_i \vartheta$, then one would not need to compute the differential consequences \eqref{eq:eg4gran:LD:difcons} of the momentum equations, and solve for the mixed velocity derivatives $\partial_t\partial_x u$, $\partial_t \partial_x v$, $\partial_t\partial_y u$, $\partial_t \partial_y v$. On the other hand, if more complex constitutive dependencies than \eqref{eq:granular2d:constit} need to be taken, then higher derivatives of the fields would appear in PDEs \eqref{eq:granular:2d:mom}, \eqref{eq:granular:2d:nu}, \eqref{eq:fluid:2d:en} and the entropy inequality \eqref{eq:entropy_ineq:eg:Gran}, and hence one would need to consider additional differential consequences, and solve for a further extended set of leading derivatives.
\end{remark}

The next step is to consider the entropy inequality \eqref{eq:entropy_ineq:eg:Gran} on the solution set. First, the entropy inequality \eqref{eq:entropy_ineq:eg:Gran} is defined and converted into a symbolic expression:
\[
\barr
\verb|Pi_S:= N(ind)*R(ind)*MaterialDer(S(Constit_Dependence2)),t)|\\
\verb|      + diff(Phi1(Constit_Dependence),x) + diff(Phi2(Constit_Dependence),y);|\\
\verb|Pi_S_Symb:=gem_analyze(Pi_S);|
\earr
\]
The entropy inequality on solutions is obtained by substituting the symbolic variable \verb|solution_set_equations_final| into \verb|Pi_S_Symb|. The numerator and denominator of the entropy inequality on solutions are computed as follows:
\[
\barr
\verb|Entropy_Inequality_Sol_Set:= combine(simplify(subs(solution_set_equations_final, Pi_S_Symb)));|\\
\verb|Entropy_Inequality_numer:=expand(numer(Entropy_Inequality_Sol_Set));|\\
\verb|Entropy_Ineqality_denom:=denom(Entropy_Inequality_Sol_Set);|
\earr
\]
In particular, the numerator is given by a rather lengthy expression, omitted here, and the denominator equals
\[
\verb|Entropy_Ineqality_denom|=\dfrac{\partial \epsilon}{\partial \theta}.
\]

The next step consists in the splitting of the entropy inequality numerator on solutions \verb|Entropy_Inequality_numer| by setting to zero all coefficients at the 82 free elements.
The latter include the independent variables, all physical fields, and all first, second and third derivatives of all fields, with the exception of the constitutive dependencies in \eqref{eq:granular2d:constit} and the prolonged set of leading derivatives
\eqref{eq:eg4gran:LD:prol}. Finally, the set of constraint equations is obtained as a symbolic variable \verb|coeffs_constraints| by the command \eqref{eg1:step:E:4}. The left-hand side of the residual entropy inequality \verb|Residual_Entropy_Ineq| is obtained the command \eqref{eg1:step:E:residual}. The final set of solution set constraints is obtained by excluding redundancies and the residual entropy inequality, by the command \eqref{eg1:step:E:5}, and stored in the variable \verb|Solution_Set_Constraints|. This yields to the total of 180 constraint equations. The residual entropy inequality and some of the constraint equations are given by lengthy expressions which are not presented here.

\subsection{Symmetry conditions. Simplification and reduction of constraint equations}

It has been previously assumed that all constitutive functions depended on all velocity gradient components $\partial_j v_i$ rather than on the symmetric strain tensor \eqref{eq:gran:constitDep:Dij}. In order to enforce the symmetry of the dependencies, the symmetry conditions are imposed additionally:
\beq\label{eq:gran:SymConds}
\dfrac{\partial \psi}{\partial (\partial_j v_i)} = \dfrac{\partial \psi}{\partial (\partial_i v_j)},
  \eeq
for every constitutive function in the set \eqref{eq:gran:constitF}:
\[
\verb|all_constit_functions:=[GEM_FREE_F[ ]];|
\]
\beq\label{eq:gran:SymConds:symb}
\barr
\verb|Symmetry_conditions:={map(q-> diff(q,Uy)=diff(q,Vx), |\\
\verb|                  all_constit_functions)[ ]} minus {0 = 0};|
\earr
\eeq
As in the previous examples, it is natural to assume the essential dependence of energy and entropy density on the temperature:
\[
\barr
\verb|Additional_assumptions:=[|\\
\verb|      diff(E(N, Nx, Ny, R, W, Wx, Wy, Ux, Uy, Vx, Vy), W)<>0,|\\
\verb|      diff(S(N, Nx, Ny, R, W, Wx, Wy, Ux, Uy, Vx, Vy), W)<>0 ];|
\earr
\]

Finally, the full set of equations and inequalities defining the physical constraint equations for the nonnegative  entropy production are given by the following set:
\beq\label{eg:gran:final:eq:ineq}
\barr
\verb|Constraints:=[Solution_Set_Constraints[], Symmetry_conditions[],|\\
\verb|              Additional_assumptions[],  Entropy_Ineqality_denom<>0 ];|
\earr
\eeq

%
%
%\[
%\barr
%\verb|general_case_rifsimp:=DEtools[rifsimp]([|\\
%\verb|     [ Solution_Set_Constraints[ ],|\\
%\verb|       Symm_conditions[ ],|\\
%\verb|       Additional_assumptions[ ] ],|\\
%\verb|     all_constit_functions, mindim=1);|\\
%\earr
%\]
%\[
%\verb|Constraint_equations_gen_case:=general_case_rifsimp[Solved];|
%\]
%

\subsection{Discussion of the results}

The solution set entropy principle applied to the granular flow model \eqref{eq:granular:eqs}, \eqref{eq:granular2d:constit} yields the set of 195 constraint equations and inequalities given by \eqref{eg:gran:final:eq:ineq}. These conditions need to be further simplified and solved, but due to complex dependencies \eqref{eq:granular2d:constit} of the constitutive functions, they do not yield a simple set of restrictions, as is the case for earlier and simpler examples considered in this paper. The detailed study of these constraint equations and inequalities requires extensive computations and analysis, and is left to future work.

In Refs.~\cite{wang_shearing_1999,wang_comparison_1999} where the classical the  M{\"u}ller-Liu procedure was used, some rather simple results were obtained, because of the following.
\begin{enumerate}
   \item Additional assumptions were used.
   \item The use of Liu identities following from the Liu's lemma resulted in setting to zero coefficients at higher derivatives separately, without taking into account differential consequences of the governing equations, which resulted in a larger number of simpler constraint equations.
\end{enumerate}
As is always the case when differential consequences need to be taken into account, computations based on Liu's lemma yield an overly-restrictive set of constraint equations. In particular, though it is not possible to present closed-form restrictions following from the solution set principle, it is possible to state that the flux relations condition
\beq\label{eq:phi:q:gran}
\Phi_i = q_i/\vartheta,\quad i=1,2,
\eeq
found to be necessary in Ref.~\cite{wang_comparison_1999} are in fact \emph{not necessary} for the nonnegative entropy production.

\end{appendix}

%------------------------------------------------------------------------------------------------------------------
\end{document}